\documentclass[a4paper,12pt,fleqn]{article}
\usepackage{amsmath,amssymb,amsthm,graphicx,subfigure,float,caption,epstopdf,tabularx,color, bm,amsfonts,epic}
\usepackage[top=1.25in, bottom=1.2in, left=1.0in, right=1.0in]{geometry}
\usepackage{appendix}
\usepackage{multirow}
\allowdisplaybreaks
\newtheorem{thm}{Theorem}[section]
\newtheorem{lem}{Lemma}[section]

\newtheorem{rmk}{Remark}[section]
\newtheorem*{prf}{Proof}
\numberwithin{equation}{section}

\usepackage{indentfirst}
\graphicspath{{Fig}}

\begin{document}

\title{A linearly implicit structure-preserving scheme for the fractional sine-Gordon equation based on the IEQ approach}

\author{Yayun Fu, \quad Wenjun Cai, \quad Yushun Wang\footnote{Correspondence author. Email:wangyushun@njnu.edu.cn.}\\{\small Jiangsu Key Laboratory for NSLSCS,}\\{\small School of Mathematical Sciences,  Nanjing Normal University,}\\{\small  Nanjing 210023, China}\\}
\date{}
\maketitle

\begin{abstract}
This paper aims to develop a linearly implicit structure-preserving numerical scheme for the space fractional sine-Gordon equation, which is based on the newly developed invariant energy quadratization method. First, we reformulate the equation as a canonical Hamiltonian system by virtue of the variational derivative of the functional with fractional Laplacian. Then, we utilize the fractional centered difference formula to discrete the equivalent system derived by the invariant energy quadratization method in space direction,
and obtain a conservative semi-discrete scheme. Subsequently, the linearly implicit structure-preserving method is applied for the resulting semi-discrete system to arrive at a fully-discrete conservative scheme. The stability, solvability and convergence in the maximum norm of the numerical scheme are given. Furthermore, a fast algorithm based on the fast Fourier transformation technique is used to reduce the computational complexity in practical computation. Finally, numerical examples are provided to confirm our theoretical analysis results.\\[2ex]


\end{abstract}

\section{Introduction}

In recent years, an increasing number of classical models that
are described by integer order partial differential equations have been formulated by the system of fractional order \cite{p1, p2, p3}. It is mainly because of that
fractional differential equations are more accurate in modelling a variety of scientific and engineering problems with long-range temporal cumulative memory effects and spatial interactions, which are frequently implemented in many fields, such as biological materials, electronic circuits, and automatic control \cite{p4, p5, p7}. The fractional sine-Gordon (FSG) equation \cite{p9,p10} is a generalization of the standard sine-Gordon (SG) equation, which represents an important dynamical model with long-range interactions \cite{p28} in nonlinear science. 
Recently, this equation has been largely studied and many significant achievements have been made in Refs. \cite{p9,p10, p14}.

In this paper, we numerically consider the following space FSG equation
\begin{align}\label{FSG:eq:1.1}
u_{tt}=-(-\Delta)^{\frac{\alpha}{2}}u-\sin(u),\ \ x\in (a,b)=\Omega,\ \ 0<t\leq T,
\end{align}
with boundary condition
\begin{align}\label{FSG:eq:1.2}
u(x,t)=0,\ \ x \in \mathbb{R}\setminus \Omega,\ \ 0\leq t \leq T，
\end{align}
and
initial conditions
\begin{align}\label{FSG:eq:1.3}
u(x,0)=\varphi(x), \ \ u_t(x,0)=\psi(x),\ \ x\in\Omega,
\end{align}
where $1< \alpha \leq 2$, $\varphi(x)$ and $\psi(x)$ are the wave modes or kinks
and their velocity, respectively. When $\alpha=2$, the FSG equation (\ref{FSG:eq:1.1}) reduces to the classical SG equation. The fractional Laplacian is defined as a pseudo-differential operator with the symbol $|\xi|^{\alpha}$ in the Fourier space \cite{p17}
\begin{align}\label{FSG:eq:1.4}
-(-\Delta)^{\frac{\alpha}{2}}u(x,t)=-\mathcal{F}^{-1}(|\xi|^{\alpha}\hat{u}(\xi,t)),
\end{align}
where $\mathcal{F}$ is the Fourier transform and $\hat{u}(\xi,t)=\mathcal{F}[u(x,t)]$. Yang \cite{p18} proposed that the fractional Laplacian is equivalent to the Riesz fractional derivative in one dimension, namely
\begin{align}\label{FSG:eq:1.5}
-(-\Delta)^{\frac{\alpha}{2}}u(x,t)=\frac{\partial^{\alpha}u(x,t)}{\partial |x|^\alpha}=-\frac{1}{2\text{cos}\frac{\alpha \pi}{2}}[_{-\infty}D^{\alpha}_x u(x,t)+ _{x}D^{\alpha}_{+\infty} u(x,t)],
\end{align}
where ${_{-\infty}}D^{\alpha}_x u(x,t)$ and $_xD^{\alpha}_{ + \infty} u(x,t)$ are the left and right side Riemann-Liouville fractional derivatives \cite{p19}, respectively, which are defined as
\begin{align}\label{FSG:eq:1.8}
_{-\infty}D^{\alpha}_x u(x,t)=\frac{1}{\Gamma (n-\alpha)}\frac{\partial^n}{\partial x^n}\int_{-\infty}^{x}(x-s)^{n-1-\alpha}u(s,t)ds,
\end{align}
\begin{align}\label{FSG:eq:1.8}
_{x}D^{\alpha}_{+\infty} u(x,t)=(-1)^{n}\frac{1}{\Gamma (n-\alpha)}\frac{\partial^n}{\partial x^n}\int_{x}^{+\infty}(s-x)^{n-1-\alpha}u(s,t)ds,
\end{align}
for $n-1< \alpha \leq n$. Some numerical schemes have been developed by scholars to approximate the Riemann-Liouville fractional derivative and Riesz fractional derivative, such as the
shifted Gr\"{u}nwald formula \cite{p20}, the difference-quadrature approach \cite{p21} and the fractional centered difference formula \cite{p22}.

Moreover, it is easy to manifest that the system \eqref{FSG:eq:1.1} with the condition \eqref{FSG:eq:1.2} possesses an energy conservation law
\begin{align}\label{FSG:eq:1.6}
 E(t)=E(0),
\end{align}
where the energy is given by
 \begin{align}\label{FSG:eq:1.7}
 E(t):=\frac{1}{2}\int_{\Omega}\Big(u_t^2+((-\Delta)^{\frac{\alpha}{4}} u)^2+2(1-\cos(u))\Big)d{x}.
 \end{align}

As is known to all, the analytical solutions of fractional differential equations contain some special functions, so it is generally difficult to obtain their explicit forms. Instead, many important numerical schemes have been discovered to solve fractional differential equations, such as finite element method \cite{p38}, finite difference method \cite{p39, p47} and spectral method \cite{p40}. In recent years, a slice of researchers have devoted great energy to discuss the computation for the FSG equation, and the readers can refer to Refs. \cite {p9, p11, p14} and references therein.

Structure-preserving algorithms are numerical methods  that can conserve one or more of the intrinsic properties of the original problem. Prior researches generally confirmed that the structure-preserving methods are more superior than the traditional methods in long time stability for numerical simulations \cite {p23, p24, p25, p26}. Thus investigating the structure-preserving algorithms for fractional differential equations exerted a tremendous fascination on researchers, which has been studied by many mathematicians and physicists \cite {p27,p28,p29,p30,p31,p32}. However, most of the available energy-preserving schemes for fractional differential equations are fully-implicit, one has to use iterations to solve a system of algebraic equations, which brings a large number of calculations in long time numerical simulation.

The invariant energy quadratization (IEQ) method developed by Yang and his collaborators, has been used to construct efficient and accurate numerical schemes for some gradient models \cite{p33,p34,p35,p36}, and the resulted schemes still retained the identical energy dissipation law. As far as we know, no previous work used this idea to investigate fractional differential equations. Moreover, few studies have focused on considering error estimates of the numerical schemes derived by the IEQ approach.
Therefore, the motivation of this paper is to extend the IEQ method to construct structure-preserving numerical scheme for fractional differential equations. To this end, taking the FSG equation as an example, we develop a linear implicit energy-preserving scheme for the equation by using the IEQ method and to estimate the error of the derived scheme.
Besides, a fast algorithm based on the fast Fourier transformation (FFT) technique is used in practical computation, which can reduce the memory requirement and the computational complexity.

The outline of this paper is as follows. In Section 2, we reformulate the FSG  equation as a canonical Hamiltonian system by virtue of the variational derivative of the functional with fractional Laplacian, and then transform the energy functional of the equation into a quadratic form of a set of new variables via a change of variables to obtain a new equivalent system in terms of the new variables. In Section 3, we use the fractional centred difference formula
to approximate the equivalent system space derivative and obtain a semi-discrete energy-preserving scheme. Then  a fully-discrete energy-preserving scheme is derived by utilizing Crank-Nicolson method to discrete the semi-discrete system in time. In Section 4, we prove that the fully-discrete scheme has a unique solution, and is  convergent with the order
of $O(h^2+\tau^2)$ in the discrete maximum norm. Numerical examples are presented in Section 5 to demonstrate the theoretical results. We draw some conclusions in Section 6.

\section{Hamiltonian formulation and IEQ method}
In this section, we derive the Hamiltonian formulation and obtain a equivalent system for the FSG equation. For each nonnegative integer $r$, let $C^{r}(\mathbb{R})$ denote the space of all the functions $\phi(x) : \mathbb{R} \rightarrow \mathbb{R}$ which have continuous derivatives up to the $r$th order, and let $L_{1}(\mathbb{R})$ represent the vector space of all the Lebesgue-integrable functions $\phi(x)$. For any function $\phi(x)\in L^2(\mathbb{R})$, we denote its Fourier transform by $\widehat{\phi}(\xi)$.

\subsection{Hamiltonian formulation and conservation law}

In this subsection, we introduce some lemmas which are extremely useful for subsequent theoretical analysis.

\begin{lem}\rm Let $s >0$, then for any real functions $p, q \in L^2({\Omega})$ with homogeneous boundary conditions, we have
\begin{align}\label{FKGS:eq:2.1}
 \big((-\Delta)^{s}p,{q}\big)=\big((-\Delta)^{\frac{s}{2}}p,{(-\Delta)^{\frac{s}{2}}q}\big)=\big(p,{(-\Delta)^{s}q}\big).
\end{align}
\end{lem}
\begin{prf}\rm  First, let us recall a useful property of Fourier transform, namely,
\begin{align}\label{FKGS:eq:2.2}
\int_{\Omega}pqdx=\int_{\Omega}\hat{p}\hat{q}d\xi.
\end{align}
Then, we can deduce
\begin{align}\label{FKGS:eq:2.3}
 \big((-\Delta)^{s}p,{v}\big)&= \big(\widehat{(-\Delta)^{s}p},\widehat{q}\big)=\big(|\xi|^{2s}\widehat{p},\widehat{q}\big)=\big(|\xi|^{s}\widehat{p},|\xi|^{s}\widehat{q}\big)\nonumber\\
&=\big(\widehat{(-\Delta)^{\frac{s}{2}}p},\widehat{(-\Delta)^{\frac{s}{2}}q}\big)=\big({(-\Delta)^{\frac{s}{2}}p},{(-\Delta)^{\frac{s}{2}}q}\big),
 \end{align}
and
\begin{align}\label{FKGS:eq:2.4}
\big((-\Delta)^{s}p,{q}\big)&= \big(\widehat{(-\Delta)^{s}p},\widehat{q}\big)=\big(|\xi|^{2s}\widehat{p},\widehat{q})=(\widehat{p},|\xi|^{2s}\widehat{q}\big)\nonumber\\
&=\big(\widehat{p},\widehat{(-\Delta)^{s}q}\big)=\big(p,{(-\Delta)^{s}q}\big).
\end{align}
\qed
\end{prf}

\begin{lem}\rm
For a functional $F[\rho]$ with the following form
\begin{align}\label{FKGS:eq:2.5}
F[\rho]=\int_{\Omega}g( \rho(\eta), (-\Delta)^{\frac{\alpha}{4}}\rho(\eta))d\eta,
\end{align}
where $g$ is a smooth function on the $\Omega$, the variational derivative of $F[\rho]$ is given as follows
\begin{align}\label{FKGS:eq:2.6}
\frac{\delta F}{\delta \rho}=\frac{\partial g}{\partial \rho}+(-\Delta)^{\frac{\alpha}{4}}\frac{\partial g}{\partial\big((-\Delta)^{\frac{\alpha}{4}}\rho\big)}.
\end{align}
\end{lem}
\begin{prf}\rm
Let $\phi(w)$ be an arbitrary function with the homogeneous boundary condition. According to the fact that the fractional Laplacian is linear, and the definition of variational derivative, we have
\begin{align}\label{FKGS:eq:2.7}
\int_{\Omega}\frac{\delta F}{\delta \rho}\phi(\eta)d\eta&=\Big[\frac{d}{d\mu}\int_{\Omega}g\big(\rho+\mu\phi,(-\Delta)^{\frac{\alpha}{4}}\rho+\mu(-\Delta)^{\frac{\alpha}{4}}\phi \big)d\eta\Big]_{\mu=0}\nonumber\\
&=\int_{\Omega}\big(\frac{\partial g}{\partial \rho} \phi+\frac{\partial g}{\partial((-\Delta)^{\frac{\alpha}{4}}\rho)}(-\Delta)^{\frac{\alpha}{4}}\phi \big)d\eta \nonumber\\
&=\int_{\Omega}\Big(\frac{\partial g}{\partial \rho} \phi+\big((-\Delta)^{\frac{\alpha}{4}}\frac{\partial g}{\partial((-\Delta)^{\frac{\alpha}{4}}\rho)}\big)\phi \Big)d\eta
\nonumber\\
&=\int_{\Omega}\Big(\frac{\partial g}{\partial \rho}+(-\Delta)^{\frac{\alpha}{4}}\frac{\partial g}{\partial((-\Delta)^{\frac{\alpha}{4}}\rho)} \Big) \phi d\eta,
\end{align}
where \eqref{FKGS:eq:2.1} was used. Based on the fact that $\phi(\eta)$ is arbitrary, by using the fundamental lemma of calculus of variations, we can obtain \eqref{FKGS:eq:2.6}.
\qed
\end{prf}

\begin{rmk}\rm
In the case of periodic boundary conditions, Lemma 2.1 and Lemma 2.2 have been proposed in Ref. \cite{p27}. In the paper,  the conclusions are given under the homogeneous boundary conditions.
\end{rmk}

Let $v=u_t$, system \eqref{FSG:eq:1.1} can be rewritten as a first-order system
\begin{align}\label{SG:eq:2.07}
u_t=v,
\end{align}
\begin{align}\label{SG:eq:2.08}
v_t=-(-\Delta)^{\frac{\alpha}{2}}u-\sin(u).
\end{align}
By taking the inner products of \eqref{SG:eq:2.07}-\eqref{SG:eq:2.08} with $v_{t}, v$, respectively, and summing them together, we can obtain system \eqref{SG:eq:2.07}-\eqref{SG:eq:2.08} has the            following energy conservation law
\begin{align*}
\frac{d}{d t}{H}=0,
\end{align*}
where the energy functional
\begin{align}\label{FKGS:eq:2.51}
{H}=\frac{1}{2}\int_{\Omega}\Big[v^{2}+((-\Delta)^{\frac{\alpha}{4}}  u)^2+2\big(1-\cos(u)\big)\Big]d{x}.
\end{align}


Based on the fractional variational  derivative formula in Lemma 2.2, we obtain the following theorem.

\begin{thm}\rm\label{2SG-lem2.2} The system \eqref{SG:eq:2.07}-\eqref{SG:eq:2.08} is an infinite-dimensional canonical Hamiltonian system
\begin{align*}
\left(\begin{array}{c}
		v_t\\
		u_t
		\end{array} \right)=J^{-1}\left(\begin{array}{c}
		\delta\mathcal{H}/\delta v\\
		\delta\mathcal{H}/\delta u
		\end{array} \right),\ \ \ J=\left(\begin{array}{cc}
		0 & 1 \\
		-1 & 0
		\end{array} \right),
\end{align*}
where the Hamiltonian function $\mathcal{H}$ is defined by
\begin{align}\label{FKGS:eq:2.50}
\mathcal{{H}}=\frac{1}{2}\int_{\Omega}\Big[v^{2}+((-\Delta)^{\frac{\alpha}{4}}  u)^2+2\big(1-\cos(u)\big)\Big]d{x}.
\end{align}

\end{thm}

%
%
%

\subsection{Invariant energy quadratization method}

We introduce an auxiliary variable $w$ in terms of $u, x, t$ with the following definition
\begin{align}\label{SG:eq:2.1}
w(u,x,t):=\sqrt{C_{0}+1-\cos(u)},
\end{align}
where $C_{0}$ is a positive constant such that $C_0+1-\cos(u) > 0$ for any $x\in \mathbb{R}$. In this paper, as an example, we take $C_0=1$. Then, energy function \eqref{FKGS:eq:2.51} of the FSG equation is transformed into a quadratic function of $w$, $u$, and $v$
\begin{align}\label{SG:eq:2.2}
{H}=\frac{1}{2}\int_{\Omega}\Big(v^{2}+((-\Delta)^{\frac{\alpha}{4}}  u)^2+2w^2\Big)d{x}.
 \end{align}
 Then, we can rewrite \eqref{SG:eq:2.07}-\eqref{SG:eq:2.08} as the following equivalent system
\begin{align}\label{SG:eq:2.3}
&u_t=v,
\end{align}
\begin{align}\label{SG:eq:2.4}
&v_t=-(-\Delta)^{\frac{\alpha}{2}}u-\frac{\sin(u)}{\sqrt{2-\cos(u)}}w,
\end{align}
\begin{align}\label{SG:eq:2.5}
&w_t=\frac{\sin(u)}{2\sqrt{2-\cos(u)}}v.
\end{align}
Taking the inner products of \eqref{SG:eq:2.3}-\eqref{SG:eq:2.5} with $v_t$, $v$ and $2w$, respectively, we derive the energy conservation law as follows
 \begin{align}\label{SG:eq:2.6}
 \frac{d}{dt}\int_{\Omega}\frac{1}{2}\Big(v^2+((-\Delta)^{\frac{\alpha}{4}}  u)^2+2w^2\Big)d{x}=0.
 \end{align}
 One can observe that the equivalent system \eqref{SG:eq:2.3}-\eqref{SG:eq:2.5} still retains a similar energy conservation law, but it is in terms of the new variable now.

\section{ Construction of the energy-preserving scheme }

In this section, we will construct an energy-preserving finite difference scheme for the FSG equation.

\subsection{Structure-preserving spatial discretization}

We choose the time step $\tau:= \frac{T}{N}$ and mesh size $h:=\frac{b-a}{M}$ with integers $N$ and $M$. Denote~$\Omega_{h}=\{x_j |\ x_j=a+jh, 0\leq j \leq M  \}$,
$\Omega_{\tau}=\{t_n|\ t_n=n\tau, 0\leq n \leq N  \}$. Let $\mathcal{V}_{h} = \{u|u=(u_1,u_2, \cdots, u_{M-1})^{T}\}$ be the space of grid functions. For a given grid function
$\mathcal{\mathring{V}}_{h} = \{u_{j}^{n}|\ u_{j}^{n}=u(x_j, t_n),  (x_j, t_n)\in \Omega_{h}\times \Omega_{\tau}\}$, we introduce some notations for any mesh function $u_{j}^{n} \in \mathcal{\mathring{V}}_{h}$ as

\begin{align*}
\delta_{t}u_{j}^{n}=\frac{u_{j}^{n+1}-u_{j}^{n}}{\tau},\ u_{j}^{n+\frac{1}{2}}=\frac{u_{j}^{n+1}+u_{j}^{n}}{2},\ \tilde{u}_{j}^{n+\frac{1}{2}}=\frac{3{ u}_{j}^{n}-{u}_{j}^{n-1}}{2}.
\end{align*}

 For any two grid functions $u$,~$v$~$\in$~$\mathcal{ V}_{h}$,~we define the discrete inner product and the associated $l^2$-norm as
 \begin{align*}
(u, v)=h\sum\limits_{j=1}^{M-1} u_{j}{v}_{j},\ \|u\|^2=(u, u),
 \end{align*}and the discrete maximum norm ($l^{\infty}$-norm) as
\begin{align*}
\|u\|_{{\infty}}=\sup\limits_{1 \leq j \leq M-1}|u_j|.
\end{align*}

Set $l_{h}^{2}$=$\{u_{j}^{n}: \|u^{n}\|^{2} < + \infty\}$. Then for $0 \leq \sigma \leq1$, $u \in l^{2}_h$, the fractional Sobolev norm $\|u\|_{H^{\sigma}}$ and the semi-norm
$|u|_{H^{\sigma}}$ can be defined as \cite{p37}
\begin{align*}
\|u\|^2_{H^{\sigma}}=h\int^{\pi}_{-\pi}(1+h^{-2\sigma}|k|^{2\sigma})|\hat{u}(k)|^2dk,\ \
|u|^2_{H^{\sigma}}=h\int^{\pi}_{-\pi}h^{-2\sigma}|k|^{2\sigma}|\hat{u}(k)|^2dk.
\end{align*}
where
\begin{align*}
\hat{u}^n(k)=\frac{1}{\sqrt{2\pi}}\sum\limits_{j\in\mathbb{Z}}u^{n}_{j}e^{-ijk}.
\end{align*}
Obviously, one can derive that
\begin{align}\label{FSG:eq:2.7}
\|u\|^{2}_{H^{\sigma}}=\|u\|^2+|u|^{2}_{H^{\sigma}}.
\end{align}

\begin{lem}\rm\cite{p22} Let $\phi(x) \in C^5(\mathbb{R})$, and all derivatives up to order five belong to $L^1(\mathbb{R})$. Then,
for $1< \alpha\leq 2$, we have
\begin{align}\label{FSG:eq:2.8}
\frac{\partial^{\alpha}\phi(x)}{\partial |x|^{\alpha}}=-\frac{1}{h^\alpha}\sum\limits_{k=-\infty}^{+\infty}c_{k}^{(\alpha)}\phi(x-kh)+O(h^2),
\end{align}
where the coefficients $c_{k}^{(\alpha)}:=\frac{(-1)^k \Gamma(\alpha+1)}{\Gamma(\frac{\alpha}{2}-k+1)\Gamma(\frac{\alpha}{2}+k+1)}$.\end{lem}
%
According to (\ref{FSG:eq:1.2})  and Lemma 3.1, we can obtain
\begin{align}\label{FSG:eq:2.10}
\frac{\partial^{\alpha}u(x,t)}{\partial |x|^{\alpha}}=-\frac{1}{h^\alpha}\sum\limits_{k=-(b-x)/h}^{-(a-x)/h}c_{k}^{(\alpha)}u(x-kh, t)+O(h^2).
\end{align}

Let $U_{j}^{n}$ and $u_{j}^{n}$ be the numerical approximation and the exact
solution of $u(x,t)$ at the points $(x_{j},t_{n})$, respectively. Then it follows from (\ref{FSG:eq:1.2}) and (\ref{FSG:eq:2.10}) that
\begin{align}\label{FSG:eq:2.11}
-(-\Delta)^{\frac{\alpha}{2}}u_{j}^{n}=-\frac{1}{h^\alpha}\sum\limits_{k=-M+j}^{j}c_{k}^{(\alpha)}u_{j-k}^{n}+O(h^2)=-\frac{1}{h^\alpha}\sum\limits_{k=1}^{M-1}c_{j-k}^{(\alpha)}u_{k}^{n}+O(h^2).
\end{align}
We introduce the notation
\begin{align}\label{FSG:eq:2.12}
\Delta_{h}^{\alpha}u_{j}^{n}=\frac{1}{h^\alpha}\sum\limits_{k=1}^{M-1}c_{j-k}^{({\alpha})}u_{k}^{n},~~~~1\leq j \leq M-1,~~0\leq n\leq N.
\end{align}

Denote matrix $\mathbf{C}$ as
\begin{align} \label{SG:eq:8.13}
\mathbf{C}=\begin{pmatrix}
       c_{0}^{(\alpha)}  & c_{-1}^{(\alpha)}& \cdots & c_{-M+2}^{(\alpha)}\\
       c_{1}^{(\alpha)}  & c_{0}^{(\alpha)}& \cdots & c_{-M+3}^{(\alpha)}\\
       \vdots & \vdots&\ddots & \vdots\\
       c_{M-2}^{(\alpha)}  & c_{M-3}^{(\alpha)}& \cdots & c_{0}^{(\alpha)}\\
       \end{pmatrix}.\end{align}
 Let $\lambda=(\lambda_1, \lambda_2, \cdots, \lambda_{M-1})$, where $\lambda_i$~($1\leq i \leq M-1$) is the eigenvalue of matrix $\mathbf{C}$, and satisfy \cite{p42}
 \begin{align}\label{SG:eq:2.13}
0<\lambda_i<2c_{0}^{(\alpha)},~~~~i=1, 2, \cdots, M-1.
\end{align}
 One can easily verify that the $\mathbf{C}$ is a real-valued symmetric positive definite Toeplitz matrix.\vskip 1mm

\begin{lem}\rm{\cite{p43}} For any two grid functions $U_1,~U_2\in \mathcal{V}_{h}$, there exists a linear operator
$\Lambda^\alpha=h^{-\frac{\alpha}{2}}\mathbf{C}^{\frac{1}{2}}$ such that
\begin{align}\label{SG:eq:2.14}
(\Delta_{h}^{\alpha}U_1,~U_2)=(\Lambda^\alpha U_1,~\Lambda^\alpha U_2),
\end{align} where the $\mathbf{C}^{\frac{1}{2}}$ is the Cholesky factor of matrix of $\mathbf{C}$, i.e., $(\mathbf{C}^{\frac{1}{2}})^2=\mathbf{C}$.\end{lem}

\begin{lem}\rm{\cite{p43}} For any grid function $U^n \in \mathcal{V}_h$, $0\leq n \leq N$, we have
\begin{align}\label{SG:eq:2.15}
(\Delta_{h}^{\alpha}U^{n+\frac{1}{2}}, \delta_{t}U^{n})=\frac{1}{2\tau}(\| \Lambda^\alpha U^{n+1} \|^{2}- \|\Lambda^\alpha U^{n}\|^{2}).
\end{align}

\subsection{Energy-preserving semi-discrete scheme}

 We now discretize system \eqref{SG:eq:2.4} in space by the fractional centred difference formula to construct a semi-discrete difference scheme, namely
\begin{align}\label{SG:eq:3.1}
&\frac{d}{dt}U_{{j}}=V_{{j}},\\ \label{SG:eq:3.2}
& \frac{d}{dt}V_{{j}}=-\Delta_{h}^{\alpha} U_{{j}}-\frac{\sin({U}_{{j}})}{\sqrt{2-\cos({U}_{{j}})}}W_{{j}},\\ \label{SG:eq:3.3}
&\frac{d}{dt}W_{{j}}=\frac{\sin({U}_{{j}})}{2\sqrt{2-\cos({U}_{{j}})}}V_{{j}}.
\end{align}

The energy conservation property of above semi-discrete scheme is given in the following theorem.

\begin{thm}\rm\label{2SG-lem2.1} The scheme \eqref{SG:eq:3.1}-\eqref{SG:eq:3.3} satisfies the semi-discrete energy conservation law
\begin{align*}
\tilde{E}(t)=\tilde{E}(0),
\end{align*}
where
\begin{align*}
\tilde{E}(t)=\frac{1}{2}\big(||{V}||^2+||\Lambda^{\alpha}{ U}||^2+2||{W}||^2).
\end{align*}
\end{thm}

\begin{prf}\rm By taking the discrete inner products of \eqref{SG:eq:3.1}-\eqref{SG:eq:3.3} with  $\frac{d}{dt}V_{{j}}$, $V_{{j}}$ and $2W_{{j}}$, respectively, one can deduce that
\begin{align}\label{SG:eq:3.4}
&(\frac{d}{dt}{U_j},\frac{d}{dt}{V_j})=\frac{1}{2}\frac{d}{dt}||{V_j}||^2,\\\label{SG:eq:3.5}
& (\frac{d}{dt}{V_j},\frac{d}{dt}{ U_j})=-\frac{1}{2}\frac{d}{dt}||\Lambda^{\alpha} { U_j}||^2-h\sum_{j=1}^{M-1}\frac{\sin({U}_{j})W_{j}}{\sqrt{2-\cos({U}_{j})}}V_{j},\\\label{SG:eq:3.6}
&\frac{d}{dt}||{W_j}||^2=h\sum_{j=1}^{M-1}\frac{\sin({U}_{j})W_{j}}{\sqrt{2-\cos({U}_{j})}}V_{j},
\end{align}
where Lemma 3.3 was used. Substituting \eqref{SG:eq:3.4} and \eqref{SG:eq:3.6} into \eqref{SG:eq:3.5}, we derive
\begin{align}\label{SG:eq:3.7}
\frac{d}{dt}(||{V}||^2+||\Lambda^{\alpha}{ U}||^2+2||{W}||^2)=0.
\end{align}
This completes the proof.
\qed
\end{prf}
\end{lem}

\subsection{A fully-discrete linear energy-preserving scheme}

Based on the discussions of the semi-discrete scheme \eqref{SG:eq:3.4}-\eqref{SG:eq:3.6} for the FSG equation, in this subsection, our goal is to establish a linear implicit fully-discrete difference scheme.

Applying the Crank-Nicolson method for \eqref{SG:eq:3.1}-\eqref{SG:eq:3.3} in time, further utilizing the extrapolation
technique, we can obtain a linear implicit scheme for the FSG equation, namely
\begin{align}\label{SG:eq:3.8}
&\delta_t^{}U_{{j}}^n=V_{{j}}^{n+\frac{1}{2}},\\\label{SG:eq:3.9}
& \delta_t^{}V_{{j}}^n=-\Delta_{h}^{\alpha}U_{{j}}^{n+\frac{1}{2}}-\frac{\sin(\tilde{U}_{{j}}^{n+\frac{1}{2}})}{\sqrt{2
-\cos(\tilde{U}_{{j}}^{n+\frac{1}{2}})}}W_{{j}}^{n+\frac{1}{2}},\\\label{SG:eq:3.10}
&\delta_t^{}W_{{j}}^n=\frac{\sin(\tilde{U}_{{j}}^{n+\frac{1}{2}})}{2\sqrt{2
-\cos(\tilde{U}_{{j}}^{n+\frac{1}{2}})}}V_{{j}}^{n+\frac{1}{2}}.
\end{align}

In addition, the first step can be obtained by some second or higher order time integrators. Here we employ the second order implicit conservative scheme as the starting step

\begin{align}\label{SG:eq:3.11}
&\delta_t^{}U_{{j}}^0=V_{{j}}^{\frac{1}{2}},\\\label{SG:eq:3.12}
& \delta_t^{}V_{{j}}^0=-\Delta_{h}^{\alpha} U_{{j}}^{\frac{1}{2}}-\frac{\sin({U}_{{j}}^{\frac{1}{2}})}{\sqrt{2-\cos({U}_{{j}}^{\frac{1}{2}})}}W_{{j}}^{\frac{1}{2}},\\\label{SG:eq:3.13}
&\delta_tW_{{j}}^0=\frac{\sin({U}_{{j}}^{\frac{1}{2}})}{2\sqrt{2-\cos({U}_{{j}}^{\frac{1}{2}})}}V_{{j}}^{\frac{1}{2}}.
\end{align}

The proposed scheme is known as IEQ-CN scheme. A result on energy conservation property of above fully-discrete scheme is presented by the
following theorem.

\begin{thm}\rm \label{SG:thm3.1} The fully-discrete scheme \eqref{SG:eq:3.8}-\eqref{SG:eq:3.13} possesses the following discrete total energy conservation law
\begin{align*}
{E}^{n+1}={E}^n,\ \ 0\leq n \leq N-1,
\end{align*}
where
 \begin{align*}
 {E}^n=\frac{1}{2}\big(||{ V}^n||^2+||\Lambda^\alpha{ U}^n||^2+2||{ W}^n||^2).
\end{align*}
\end{thm}
\begin{prf}\rm By taking the discrete inner products of \eqref{SG:eq:3.8}, \eqref{SG:eq:3.9} and \eqref{SG:eq:3.10} with $\delta_{t} V_{{j}}^n$, $V_{{j}}^{n+\frac{1}{2}}$ and $W_{{j}}^{n+\frac{1}{2}}$, respectively, one gets
\begin{align}\label{SG:eq:3.14}
&(\delta_t{ U_j}^n,\delta_t{V_j}^n) =\frac{1}{2}\delta_t||{ V_j}^n||^2,\\\label{SG:eq:3.15}
& (\delta_t{ V_j}^n,\delta_t{ U_j}^n)=-\frac{1}{2}\delta_t||\Lambda^\alpha { U_j}^n||^2-h\sum_{j=1}^{M-1}\frac{\sin(\tilde{U}_{j}^{n+\frac{1}{2}})W_{j}^{n+\frac{1}{2}}}{\sqrt{2-\cos(\tilde{U}_{j}^{n+\frac{1}{2}})}}V_{j}^{n+\frac{1}{2}},\\\label{SG:eq:3.16}
&\frac{1}{2}\delta_t||{ W_j}^n||^2=\frac{1}{2}h\sum_{j=1}^{M-1}\frac{\sin(\tilde{U}_{j}^{n+\frac{1}{2}})W_{j}^{n+\frac{1}{2}}}{\sqrt{2-\cos(\tilde{U}_{j}^{n+\frac{1}{2}})}}V_{j}^{n+\frac{1}{2}}.
\end{align}
Substituting \eqref{SG:eq:3.14} and \eqref{SG:eq:3.16} into \eqref{SG:eq:3.15}, we have
\begin{align}\label{SG:eq:3.17}
||{ V}^{n+1}||^2-||{ V}^n||^2+||\Lambda^\alpha{ U}^{n+1}||^2-||\Lambda^\alpha{ U}^n||^2+2||{ W}^{n+1}||^2-2||{ W}^n||^2=0,
\end{align}
which implies that
\begin{align}\label{SG:eq:3.18}
{E}^{n+1}={E}^n, \ \ 1\le n\le N-1.
\end{align}

Similarly, we compute the discrete inner product of  \eqref{SG:eq:3.11}, \eqref{SG:eq:3.12} and \eqref{SG:eq:3.13} with $\delta_tV_{{j}}^0$, $V_{{j}}^{\frac{1}{2}}$ and $W_{{j}}^{\frac{1}{2}}$, respectively. We can prove
\begin{align}\label{SG:eq:3.19}
{E}^1={E}^{0}.
\end{align}
This, together with \eqref{SG:eq:3.18}, we can obtain
\begin{align*}
{E}^{n+1}={E}^n, \ \ 1\le n\le N-1.
\end{align*}
We complete the proof.
\qed
\end{prf}

\section{Numerical analysis }

In this section, we discuss the stability, solvability and convergence of the IEQ-CN scheme.

\subsection{Stability and solvability}

 According to the Theorem \ref{SG:thm3.1} and the discrete conservation law of energy, we first present the following stability result of the IEQ-CN scheme.

\begin{thm} \rm If $\varphi(x) \in H_{0}^{1}(\Omega), \psi(x) \in L^{2}(\Omega),$  there exits some positive constants $C$ such that
\begin{align}\label{SG:eq:3.20}
||{ V}^n||\le C,\ \ ||\Lambda^\alpha{ U}^n||\le C, \ \ ||W^n||\le C, \ \ n=0, 1, 2, \cdots,
\end{align}
which implies that the IEQ-CN scheme is stable.
\end{thm}

Here and subsequent
theoretical analysis, $C$ is a general positive constant which is independent of $\tau$ and $h$. Note that $C
$ may vary in
different circumstances.

We present the uniquely solvable property of the IEQ-CN  scheme
in the following theorem.

\begin{thm}\rm \label{2SG:thm3.2} The fully-discrete difference scheme \eqref{SG:eq:3.8}-\eqref{SG:eq:3.10} is uniquely solvable.\end{thm}
\begin{prf}\rm
One can observe that
\begin{align}\label{SG:eq:4.1}
V_{{j}}^{n+\frac{1}{2}}=\frac{U_{{j}}^{n+1}-U_{{j}}^{n}}{\tau}=\frac{2U_{{j}}^{n+\frac{1}{2}}-2U_{{j}}^{n}}{\tau},
\end{align}
\begin{align}\label{SG:eq:4.2}
\delta_{t}W_{{j}}^{n}=\frac{W_{{j}}^{n+1}-W_{{j}}^{n}}{\tau}=\frac{2W_{{j}}^{n+\frac{1}{2}}-2W_{{j}}^{n}}{\tau}.
\end{align}
From \eqref{SG:eq:3.10} and \eqref{SG:eq:4.2}, we have
\begin{align}\label{SG:eq:4.3}
 W_{{j}}^{n+\frac{1}{2}}&=W_{{j}}^n
+\frac{\sin(\tilde{U}_{{j}}^{n+\frac{1}{2}})}{2\sqrt{2-\cos(\tilde{U}_{{j}}^{n+\frac{1}{2}})}}\Big(U_{{j}}^{n+\frac{1}{2}}-U_{{j}}^{n}\Big).
\end{align}
Together with \eqref{SG:eq:4.1} and \eqref{SG:eq:4.3}, system \eqref{SG:eq:3.9} can be written as
\begin{align}\label{SG:eq:4.4}
U_{{j}}^{n+\frac{1}{2}}=&U_{{j}}^{n}+\frac{\tau}{2}V_{{j}}^{n}-\frac{\tau^2}{4}\Delta_{h}^{\alpha} U_{{j}}^{n+\frac{1}{2}}\nonumber\\
&-\frac{\tau^2}{4}\frac{\sin(\tilde{U}_{{j}}^{n+\frac{1}{2}})}{\sqrt{2
-\cos(\tilde{U}_{{j}}^{n+\frac{1}{2}})}}\Big[W_{{j}}^n
+\frac{\sin(\tilde{U}_{{j}}^{n+\frac{1}{2}})}{2\sqrt{2-\cos(\tilde{U}_{{j}}^{n+\frac{1}{2}})}}\Big(U_{{j}}^{n+\frac{1}{2}}
-U_{{j}}^{n}\Big)\Big].
\end{align}

We denote matrix $\mathcal{A}$ as follows
\begin{align*}
\mathcal{A}=\left(\begin{array}{ccccc}
              {I} & -\frac{\tau}{2}{I} &0 \\
              -\frac{\tau}{2}{\Delta}_{h}^{\alpha}& {I} &\frac{\tau}{2}\text{diag}({\mathcal{B}}(\tilde{U}^{n+\frac{1}{2}}))\\
              0 &-\frac{\tau}{2}\text{diag}({\mathcal{B}}(\tilde{U}^{n+\frac{1}{2}})) & 2{I}
             \end{array}
\right),
\end{align*}
where ${I} $ is a unit matrix. Then, the scheme \eqref{SG:eq:3.8}-\eqref{SG:eq:3.10} can be reformed as the following linear equations
\begin{align} \label{SG:eq:4.80}
\mathcal{A}{ Z}^{n+\frac{1}{2}}={\mathcal{B}}^n,\ {Z}=({U}^T, {V}^T, {W}^T)^T,\ {\mathcal{B}}^n=(({U}^n)^T, ({ V}^n)^T, 2({W}^n)^T)^T,
\end{align}
where
$${\mathcal{B}}(\tilde{U}^{n+\frac{1}{2}})=\Big({\mathcal{B}}(\tilde{U}_{0}^{n+\frac{1}{2}}),\cdots,{\mathcal{B}}(\tilde{U}_{M-1}^{n+\frac{1}{2}})\Big)^T,\ \mathcal{B}(\tilde{U}_{{j}}^{n+\frac{1}{2}})=\frac{\sin(\tilde{U}_{{j}}^{n+\frac{1}{2}})}{\sqrt{2-\cos(\tilde{U}_{{j}}^{n+\frac{1}{2}})}}.$$
It is clear to see that
\begin{align*}
\mathcal{A}&=\left(\begin{array}{ccccc}
              0 & \frac{\tau}{4}{(\Delta_{h}^{\alpha}}-I) &0 \\
              -\frac{\tau}{4}{(\Delta_{h}^{\alpha}}-I)& 0 &\frac{\tau}{2}\text{diag}({\mathcal{B}}(\tilde{U}^{n+\frac{1}{2}}))\\
              0 &-\frac{\tau}{2}\text{diag}({\mathcal{B}}(\tilde{U}^{n+\frac{1}{2}})) & 0
             \end{array}
\right)\\&~~~~+\left(\begin{array}{ccccc}
              {I}& -\frac{\tau}{4}{(I+\Delta_{h}^{\alpha}})&0\\
              -\frac{\tau}{4}{(I+\Delta_{h}^{\alpha}})& {I} &0\\
              0 &0 & 2{I}
             \end{array}
\right)\\
&=:\mathcal{A}_1+\mathcal{A}_2.
\end{align*}
When $\mathcal{ A}{ x}=0$, based on the anti-symmetry of $\mathcal{A}_1$, we have
\begin{align}\label{SG:eq:4.5}
0={x}^T\mathcal{ A}{ x}={x}^T\big(\mathcal{ A}_1+\mathcal{A}_2\big){x}={ x}^T\mathcal{ A}_2{ x}.
\end{align}

Let $\mathcal{C}$ be the second-order leading principle minor of the matrix $\mathcal{ A}$, and we can get
\begin{align*}
 \mathcal{C}=\left(\begin{array}{ccccc}
              {I}& -\frac{\tau}{4}{(I+\Delta_{h}^{\alpha}})\\
              -\frac{\tau}{4}{(I+\Delta_{h}^{\alpha}})& {I} \\
              \end{array}
\right).
\end{align*}
When $\tau$ is a sufficiently small positive constant, based on \eqref{SG:eq:2.13}, we deduce that
\begin{align}\label{SG:eq:4.6}
|\mathcal{C}|=|I-\frac{\tau^2}{16}(I+\Delta_{h}^{\alpha})^2|>0,
\end{align}
which leads that the matrix $\mathcal{A}_{2}$ is symmetric positive definite.

The above discussions  demonstrate that the solutions of $\mathcal{ A}_{2}{ x}=0$ are the solutions of $x^{T}\mathcal{ A}_2{ x}=0$.

Noting that the $\mathcal{ A}_2{ x}=0$ has only zero solution, which implies that $\mathcal{ A}_2$ is invertible. From \eqref{SG:eq:4.5}, we can obtain that $\mathcal{ A}{ x}=0$ has only zero solution. By noting ${Z}^{n+1}=2{Z}^{n+\frac{1}{2}}-{Z}^n$, we finish the proof.
\qed
\end{prf}

\subsection{Convergence analysis}

In this subsection, we will establish an optimal priori estimate for the proposed scheme in discrete $l^{\infty}$-norm.
Some basic lemmas for subsequent theoretical analysis are given as follows.

\begin{lem}\rm \label{SG:lem4.1}
 (Gronwall Inequality \cite{p44}).
 Assume that the discrete function $\{s^{n}|n=0,1, \cdots, N: N\tau=T\}$ satisfy the following recurrence formula
\begin{align*}
s^n-s^{n-1}\leq \tilde{A} \tau s^n+\tilde{B} \tau s^{n-1}+\tilde{C_n}\tau,\end{align*}
where $\tilde{A}, \tilde{B}, \tilde{C_n} (n=1, 2, \cdots, N)$ are nonnegative constants. Then
 \begin{align*}\max\limits_{1 \leq n \leq N}|s^n|\leq (s^0+\tau\sum\limits_{k=1}^{N}\tilde{C_k})e^{2(\tilde{A}+\tilde{B})},
 \end{align*} where $\tau$ is a sufficiently small positive constant with $(\tilde{A}+\tilde{B})\tau \leq \frac{N-1}{2N}, N>1$.
\end{lem}

\begin{lem}\rm \label{2SG:lem4.2}(Discrete Sobolev inequality \cite{p37}). Note that $\|u^n\|_{H^{\sigma}}$ $\leq$ $\|u^n\|_{H^{\rho}}$ for $0 \leq \sigma \leq \rho \leq 1$. Then for every $\frac{1}{2}<\sigma\leq1$, there exists a constant $C=C(\sigma)> 0$ independent of $h > 0$, such that
\begin{align*}
\|u^n\|_{\infty} \leq C\|u^n\|_{H^{\sigma}},
\end{align*}
for all $u^n \in l_{h}^{2}$.\end{lem}

\begin{lem}\rm \label{2SG:lem4.3}(Uniform norm equivalence \cite{p30}). For every $1 <\alpha\leq 2$, we have
\begin{align*}
|\frac{2}{\pi}|^{\alpha}|u^n|^2_{H^{\alpha/2}}\leq h\sum\limits_{j=-\infty}^{\infty}(\triangle_{h}^{\alpha}u_{j}^{n}{u_{j}^{n}})\leq|u^n|^2_{H^{\alpha/2}},
\end{align*}
and
\begin{align*}
|\frac{2}{\pi}|^{\alpha}|u^n|_{H^{\alpha/2}}|v^n|_{H^{\alpha/2}}\leq h\sum\limits_{j=-\infty}^{\infty}(\triangle_{h}^{\alpha}u_{j}^{n}{v_{j}^{n}})\leq|u^n|_{H^{\alpha/2}}|v^n|_{H^{\alpha/2}}.
\end{align*}
\end{lem}

\begin{lem}\rm \label{2SG:lem4.4} Let $\mathcal{B}(x)=\frac{\sin(x)}{\sqrt{2-\cos(x)}},$  for any $ x\in\mathbb{R}$, one can obtain that there exist constants $C$, such that

$$
|\mathcal{B}(x)| \leq 1,\ \ \left|\mathcal{B}^{'}(x)\right|=\left|\frac{\cos (x)}{\sqrt{2-\cos x}}-\frac{\sin ^{2} (x)}{2\big(2-\cos (x)\big)^{\frac{3}{2}}}\right| \leq \frac{3}{2}
$$
and
$$
\left|\mathcal{B}^{''}(x)\right|=\left|-\mathcal{B}(x)-\frac{3 \sin 2 (x)}{4\big(2-\cos (x)\big)^{\frac{3}{2}}}+\frac{3 \sin ^{3} (x)}{4\big(2-\cos (x)\big)^{\frac{5}{2}}}\right| \leq \frac{5}{2}.
$$

\end{lem}

Then the convergence of the fully-discrete scheme is given by the following theorem.
\begin{thm}\rm \label{2SG:thm4.2} Let $u(x,t)\in C^{4}\big([0,T];C^{5}(\mathbb{R})\cap L^{1}(\mathbb{R})\big)$ be the exact solution of the original problem \eqref{FSG:eq:1.1}, and ${U}^{n}$ be the numerical solutions of the IEQ-CN scheme \eqref{SG:eq:3.8}-\eqref{SG:eq:3.10} at time level $n$. Then, as $\tau$ is sufficiently small, we have
\begin{align*}
&||{ u}^{n}-{U}^{n}||_{\infty}\leq C(h^2+\tau^{2}).
\end{align*}
\end{thm}

\begin{prf} \rm

Let $(U_{j}^{n}, V_{j}^{n}, W_{j}^{n} )$ denote the numerical approximation to the exact solution $(u, v, w)$ at the point $(x_j,t_n)$. First, we can easily prove that
\begin{align}\label{SG:eq:4.35}
||{ u}^{1}-{U}^{1}||_{\infty}\leq C(h^2+\tau^{2}).
\end{align}

Define the truncation errors of the scheme \eqref{SG:eq:3.8}-\eqref{SG:eq:3.10} as follows
\begin{align}\label{SG:eq:4.7}
&\delta_tu_{{j}}^n=v_{{j}}^{n+\frac{1}{2}}+(q_1)_{{j}}^n,
\end{align}
\begin{align}\label{SG:eq:4.8}
& \delta_tv_{{j}}^n=-\Delta_{h}^{\alpha} u_{{j}}^{n+\frac{1}{2}}-\frac{\sin(\tilde{u}_{{j}}^{n+\frac{1}{2}})}{\sqrt{2-\cos(\tilde{u}_{{j}}^{n+\frac{1}{2}})}}w_{{j}}^{n+\frac{1}{2}}+(q_2)_{{j}}^n,
\end{align}
\begin{align}\label{SG:eq:4.9}
&\delta_tw_{{j}}^n=
\frac{\sin(\tilde{u}_{{j}}^{n+\frac{1}{2}})}{2\sqrt{2-\cos(\tilde{u}_{{j}}^{n+\frac{1}{2}})}}v_{{j}}^{n+\frac{1}{2}}+(q_3)_{{j}}^n,\ 1\le n\le N-1.
\end{align}
By the Taylor expansion and \eqref{FSG:eq:2.11}, we can obtain
\begin{align}\label{SG:eq:4.10}
||{q}_1^n||\le C(h^2+\tau^2),\end{align}
\begin{align}\label{SG:eq:4.11}
||{ q}_2^n||\le C(h^2+\tau^2),
\end{align}
\begin{align}\label{SG:eq:4.12}
||{q}_3^n||\le C(h^2+\tau^2),\ \ \ 1\le n\le N-1.
\end{align}

Let
\begin{align*}
(\varepsilon_1)_{{j}}^n=u_{{j}}^n-U_{{j}}^n,\ (\varepsilon_2)_{{j}}^n=v_{{j}}^n-V_{{j}}^n,\ (\varepsilon_3)_{{j}}^n=w_{{j}}^n-W_{{j}}^n.
\end{align*}
Subtracting \eqref{SG:eq:4.7}-\eqref{SG:eq:4.9} from \eqref{SG:eq:3.8}-\eqref{SG:eq:3.10}, respectively, yields that
\begin{align}\label{SG:eq:4.13}
&\delta_t(\varepsilon_1)_{{j}}^n=(\varepsilon_2)_{{j}}^{n+\frac{1}{2}}+(q_1)_{{j}}^n,
\end{align}
\begin{align}\label{SG:eq:4.14}
& \delta_t(\varepsilon_2)_{{j}}^n=-\Delta^{\alpha}_h (\varepsilon_1)_{{j}}^{n+\frac{1}{2}}-\mathcal{B}(\tilde{u}_{{j}}^{n+\frac{1}{2}})w_{{j}}^{n+\frac{1}{2}}
+\mathcal{B}(\tilde{U}_{{j}}^{n+\frac{1}{2}})W_{{j}}^{n+\frac{1}{2}}+(q_2)_{{j}}^n,
\end{align}
\begin{align}\label{SG:eq:4.15}
&\delta_t(\varepsilon_3)_{{j}}^n=
\frac{1}{2}\mathcal{B}(\tilde{u}_{{j}}^{n+\frac{1}{2}})v_{{j}}^{n+\frac{1}{2}}
-\frac{1}{2}\mathcal{B}(\tilde{U}_{{j}}^{n+\frac{1}{2}})V_{{j}}^{n+\frac{1}{2}}+(q_3)_{{j}}^n.
\end{align}
Taking discrete inner product of \eqref{SG:eq:4.13}-\eqref{SG:eq:4.15} with ${\varepsilon}_1^{n+\frac{1}{2}}$, ${\varepsilon}_2^{n+\frac{1}{2}}$ and
${\varepsilon}_3^{n+\frac{1}{2}}$, respectively, we get

\begin{align}\label{SG:eq:4.16}
\frac{1}{2}\delta_t||{\varepsilon}_1^n||^2=({\varepsilon}_2^{n+\frac{1}{2}},{\varepsilon}_1^{n+\frac{1}{2}})+({ q}_1^n,{\varepsilon}_1^{n+\frac{1}{2}}),
\end{align}

\begin{align}\label{SG:eq:4.17}
\frac{1}{2}\delta_t||{ \varepsilon}_2^n||^2+\frac{1}{2}\delta_t||\Lambda^{\alpha}{ \varepsilon}_1^n||^2=&({ q}_2^n,{\varepsilon}_2^{n+\frac{1}{2}})+({\Delta_{h}^{\alpha} {\varepsilon}}_1^{n+\frac{1}{2}},{ q}_1^n)-\Big(\mathcal{B}({ \tilde{U}}^{n+\frac{1}{2}})\cdot { \varepsilon}_3^{n+\frac{1}{2}},{\varepsilon}_2^{n+\frac{1}{2}}\Big)\nonumber\\&-\Big( (\mathcal{B}(\tilde{ u}^{n+\frac{1}{2}})-\mathcal{B}(\tilde{ U}^{n+\frac{1}{2}}))\cdot{ w}^{n+\frac{1}{2}},{ \varepsilon}_2^{n+\frac{1}{2}}\Big)
\end{align}

\begin{align}\label{SG:eq:4.18}
\frac{1}{2}\delta_t{}||{ \varepsilon}_3^n||^2=&\frac{1}{2}\Big((\mathcal{B}(\tilde{u}^{n+\frac{1}{2}})-\mathcal{B}(\tilde{ U}^{n+\frac{1}{2}}))\cdot{ v}^{n+\frac{1}{2}},{ \varepsilon}_3^{n+\frac{1}{2}}\Big)+({ q}_3^n,{\varepsilon}_3^{n+\frac{1}{2}})\nonumber\\&+\frac{1}{2}\Big(\mathcal{B}({ U}^{n+\frac{1}{2}})\cdot {\varepsilon}_2^{n+\frac{1}{2}},{ \varepsilon}_3^{n+\frac{1}{2}}\Big),
\end{align}
where `$\cdot$' means the point multiplication between vectors, i.e., $u\cdot v=(u_{0}v_{0}, \cdots, u_{M}v_{M} )^{T}$.

Based on Lemma 4.4, one immediately obtains that
\begin{align}\label{SG:eq4.19}
||\Big(\mathcal{B}(\tilde{ u}^{n+\frac{1}{2}})-\mathcal{B}(\tilde{U}^{n+\frac{1}{2}})\Big)\cdot{w}^{n+\frac{1}{2}}||\le C||\tilde{ \varepsilon}_1^{n+\frac{1}{2}}||,
\end{align}
\begin{align}\label{SG:eq:4.20}
||\Big(\mathcal{B}(\tilde{u}^{n+\frac{1}{2}})-\mathcal{B}(\tilde{U}^{n+\frac{1}{2}})\Big)\cdot{v}^{n+\frac{1}{2}}||\le C||\tilde{ \varepsilon}_1^{n+\frac{1}{2}}||,
\end{align}
\begin{align}\label{SG:eq:4.21}
||\mathcal{B}(\tilde{ U}^{n+\frac{1}{2}})\cdot {\varepsilon}_3^{n+\frac{1}{2}}||\le C||{ \varepsilon}_3^{n+\frac{1}{2}}||,
\end{align}
\begin{align}\label{SG:eq:4.22}
||\mathcal{B}(\tilde{U}^{n+\frac{1}{2}})\cdot { \varepsilon}_2^{n+\frac{1}{2}}||\le C||{\varepsilon}_2^{n+\frac{1}{2}}||.
\end{align}
The above discussions indicate that
\begin{align}\label{SG:eq:4.23}
\delta_t||{\varepsilon}_1^n||^2&\le C\big(||{ \varepsilon}_1^{n}||^2+||{\varepsilon}_1^{n+1}||^2+||{ \varepsilon}_2^{n}||^2+||{\varepsilon}_2^{n+1}||^2
\big)+C(h^2+\tau^2)^2,
\end{align}
\begin{align}\label{SG:eq:4.24}
\delta_t||{ \varepsilon}_2^n||^2+\delta_t||\Lambda^{\alpha}{ \varepsilon}_1^n||^2\le C&\big( ||{\varepsilon}_1^{n-1}||^2+||{ \varepsilon}_1^{n}||^2+||\Lambda^{\alpha}{\varepsilon}_1^{n}||^2+||\Lambda^{\alpha}{\varepsilon}_1^{n+1}||^2\nonumber\\
&+||{\varepsilon}_2^{n}||^2+||{\varepsilon}_2^{n+1}||^2+||{\varepsilon}_3^{n}||^2+||{ \varepsilon}_3^{n+1}||^2\big)\nonumber\\&+C(h^2+\tau^2)^2,
\end{align}
\begin{align}\label{SG:eq:4.25}
\delta_t{}||{\varepsilon}_3^n||^2\le C&\big(||{\varepsilon}_1^{n-1}||^2+||{\varepsilon}_1^{n}||^2+||{\varepsilon}_2^{n}||^2
+||{\varepsilon}_2^{n+1}||^2+||{\varepsilon}_3^{n}||^2+||{\varepsilon}_3^{n+1}||^2\big)\nonumber\\&+C(h^2+\tau^2)^2.
\end{align}
According to \eqref{SG:eq:4.23}-\eqref{SG:eq:4.25}, one gets
\begin{align}\label{SG:eq4.26}
\delta_t\big(||\Lambda^{\alpha}{ \varepsilon}_1^n||^2+||{ \varepsilon}_1^n||^2+||{\varepsilon}_2^n||^2+&||{\varepsilon}_3^n||^2\big)\le C\big( ||{\varepsilon}_1^{n-1}||^2+||{ \varepsilon}_1^{n}||^2+||\Lambda^{\alpha}{\varepsilon}_1^{n}||^2\nonumber\\
&+||\Lambda^{\alpha}{\varepsilon}_1^{n+1}||^2+||{\varepsilon}_2^{n}||^2+||{\varepsilon}_2^{n+1}||^2+||{\varepsilon}_3^{n}||^2\nonumber\\&
+||{ \varepsilon}_3^{n+1}||^2\big)+C(h^2+\tau^2)^2.
\end{align}
Together with Lemma 4.3, we can deduce that
\begin{align}\label{SG:eq4.27}
\mathcal{E}^{n+1}-\mathcal{E}^{n}
 \le &C\tau\big( ||{\varepsilon}_1^{n-1}||^2+||{ \varepsilon}_1^{n}||^2+||{\varepsilon}_1^{n}||^2+||{\varepsilon}_1^{n+1}||^2+||{\varepsilon}_2^{n}||^2\nonumber\\
&+||{\varepsilon}_2^{n+1}||^2+||{\varepsilon}_3^{n}||^2+||{ \varepsilon}_3^{n+1}||^2\big)+C\tau(h^2+\tau^2)^2,
\end{align}
where the $\mathcal{E}^{n}$ ($n\geq2$) is defined by
\begin{align}\label{SG:eq4.28}
\mathcal{E}^{n}:=|{\varepsilon}_1^{n}|_{H^{\alpha/2}}^{2}+||{\varepsilon}_1^{n}||^2+||{\varepsilon}_2^{n}||^2+||{\varepsilon}_3^{n}||^2.
\end{align}
According to Lemma 4.1, one gets
\begin{align}\label{SG:eq:4.29}
\mathcal{E}^{n}-\mathcal{E}^{n-1}&\le C\tau (\mathcal{E}^{n}+\mathcal{E}^{n-1})+C\tau||{\varepsilon}_1^{n-2}||^2+C\tau(h^2+\tau^2)^2.
\end{align}
Summing up for the superscript $n$ from 2 to $m$ and then replacing $m$ by $n$, we can deduce from Lemma 4.1 that
\begin{align}\label{SG:eq:4.30}
\mathcal{E}^{n}&\le C\tau \sum_{k=1}^n\mathcal{E}^{k}+C\tau||{\varepsilon}_1^0||^2+CT(h^2+\tau^2)^2\nonumber\\&= C\tau \sum_{k=1}^n\mathcal{E}^{k}+CT(h^2+\tau^2)^2,
\end{align}
where ${\varepsilon}_1^0={0}$ is used. Applying Lemma \ref{SG:lem4.1}, we can arrive at
\begin{align}\label{SG:eq:4.31}
 |{\varepsilon}_{1}^n|_{H^{\alpha/2}}^2+||{\varepsilon}_{1}^n||^2+||{\varepsilon}_2^n||^2+||{\varepsilon}_3^n||^2\le C e^{2CT}(h^2+\tau^2)^2.
 \end{align}
It is observe that
 \begin{align}\label{SG:eq:4.32}
 |{\varepsilon}_{1}^n|_{H^{\alpha/2}}^2+||{\varepsilon}_{1}^n||+||{\varepsilon}_2^n||+||{\varepsilon}_3^n||\le C(h^2+\tau^2),
 \end{align}
 which further implies that
\begin{align}\label{SG:eq:4.33}
||{\varepsilon}_{1}^n||+|{\varepsilon}_1^n|_{H^{\alpha/2}}\le C( h^2+\tau^2).
\end{align}
Together with \eqref{FSG:eq:2.7} and Lemma 4.2, it is easy to verify straightforwardly that
\begin{align}\label{SG:eq:4.34}
||{\varepsilon}_1^n||_{\infty}=||{ u}^{n}-{U}^{n}||_{\infty}\le C( h^2+\tau^2), ~~2\leq n \leq N.
\end{align}
This completes the proof.
\qed
\end{prf}

\section{Numerical examples}
In this section,  numerical examples of the IEQ-CN scheme \eqref{SG:eq:3.8}-\eqref{SG:eq:3.10} are presented to illustrate the previous theoretical results. In the simulations, we noticed that the matrix \textbf{C} defined in \eqref{SG:eq:8.13} is a Toeplitz matrix. To reduce the memory requirement and the computational complexity in practical computation, a fast algorithm based on the FFT technique  \cite{p48} is used to solve the linear system \eqref{SG:eq:4.80}.

Let $E(h,\tau)$ be the error function between the numerical solution $U(x,t)$ and the analytical solution $u(x,t)$ at the point $(h,\tau)$, which is defined as
\begin{align}\label{FSG:eq:5.1}
E(h,\tau)=\|u(h,\tau)-U(h,\tau)\|_{\infty}.
\end{align}
When $1< \alpha < 2$, the
exact solution is not given. To obtain numerical errors, we use the error function defined as follows
\begin{align}\label{FSG:eq:5.2}
E(h,\tau)=\|U_{M}^{N}(h,\tau)-U_{2M}^{2N}(h/2,\tau/2)\|_{\infty},
\end{align}
where $h, \tau$ are mesh and time steps, respectively.
For calculating the convergence order, we use the formula
\begin{align*}
{p}=\text{log}_{2}{(E(h,\tau)/E(h/2,\tau/2))}.
\end{align*}
The relative energy error
is defined as
\begin{align*}
RE^{n}=|(E^{n}-E^{0})/E^{0}|,
\end{align*}
where $E^{n}$ denotes the energy at $t=n\tau$.

\textbf{Example 5.1.} We study system (\ref{FSG:eq:1.1}) with different fractional order $\alpha$. The initial conditions are chosen as
\begin{align}
&\varphi(x)=0,\\
&\psi(x)=\frac{4}{\omega}\text{sech}\Big(\frac{x}{\omega}\Big),\ x\in \Omega.
\end{align}
In our computation,  we take the Dirichlet boundary condition as
\begin{align*}
u(x,t)=0,\ \ x \in \mathbb{R}\setminus \Omega,\ \ 0\leq t \leq T.
\end{align*}

When $\alpha=2$, system (\ref{FSG:eq:1.1}) reduces to the standard nonlinear SG equation with the analytical solution given by
\begin{align}
u(x,t)=4\tan^{-1}\Big[\phi(t;\omega)\text{sech}\Big(\frac{x}{\omega}\Big)\Big],
\end{align}
where
\begin{align}
\phi(t;\omega)=\left\{\begin{array}{lll}
              &\frac{\sin(\omega^{-1}\sqrt{\omega^2-1}t)}{\sqrt{\omega^2-1}},\ &\text{if}\ \omega>1,\\
              \\
              &~~~~~~~~~~~~~~~~t,\ &\text{if}\ \omega=1,\\
              \\
              &\frac{\sinh(\omega^{-1}\sqrt{1-\omega^2}t)}{\sqrt{1-\omega^2}},\ &\text{if}\ 0<\omega<1.\\
             \end{array}
\right.
\end{align}

First, we test the convergence rates and the efficiency of the IEQ-CN  scheme. In our computation, we set the space interval $\Omega=(-20,20)$. Without loss of generality, we take $\omega=1.1$ and test the convergence orders of the IEQ-CN scheme for different $\alpha$. Table. 1 shows the errors and convergence orders, which
indicates that our scheme is of second-order accuracy in both space and time, which confirms the theoretical analysis.
The motivation of our work is to develop a more efficient structure-preserving scheme, thus, it is valuable to compare our new scheme with some existing schemes in computing efficiency, as follows:

\begin{itemize}
 \item D-IEQ: The direct algorithm to  solve the linear system \eqref{SG:eq:4.80} by using the preconditioned conjugate gradients method \cite{p49}.

  \item F-IEQ: The fast algorithm based on the FFT technique is applied for solving the linear system \eqref{SG:eq:4.80}.

  \item I-FDS: The implicit energy-preserving difference schemes stated in Refs. \cite{p10,p11,p28} for the one-dimensional sine-Gordon equation.
\end{itemize}
We use the standard fixed-point iteration for the fully-implicit schemes and set $10^{-14}$ as the error tolerance for all the problems. The consumed CPU time of different methods solving the FSG equation are displayed in Fig.1.
Numerical experiments show that the cost of I-FDS is most expensive while the one of F-IEQ is cheapest.
Therefore, it is preferable to construct linear implicit schemes through the IEQ approach and develop corresponding fast algorithms for large scale simulations, keeping the system energy being preserved as well.

\begin{figure}[H]
\centering\begin{minipage}[t]{70mm}
\includegraphics[width=70mm]{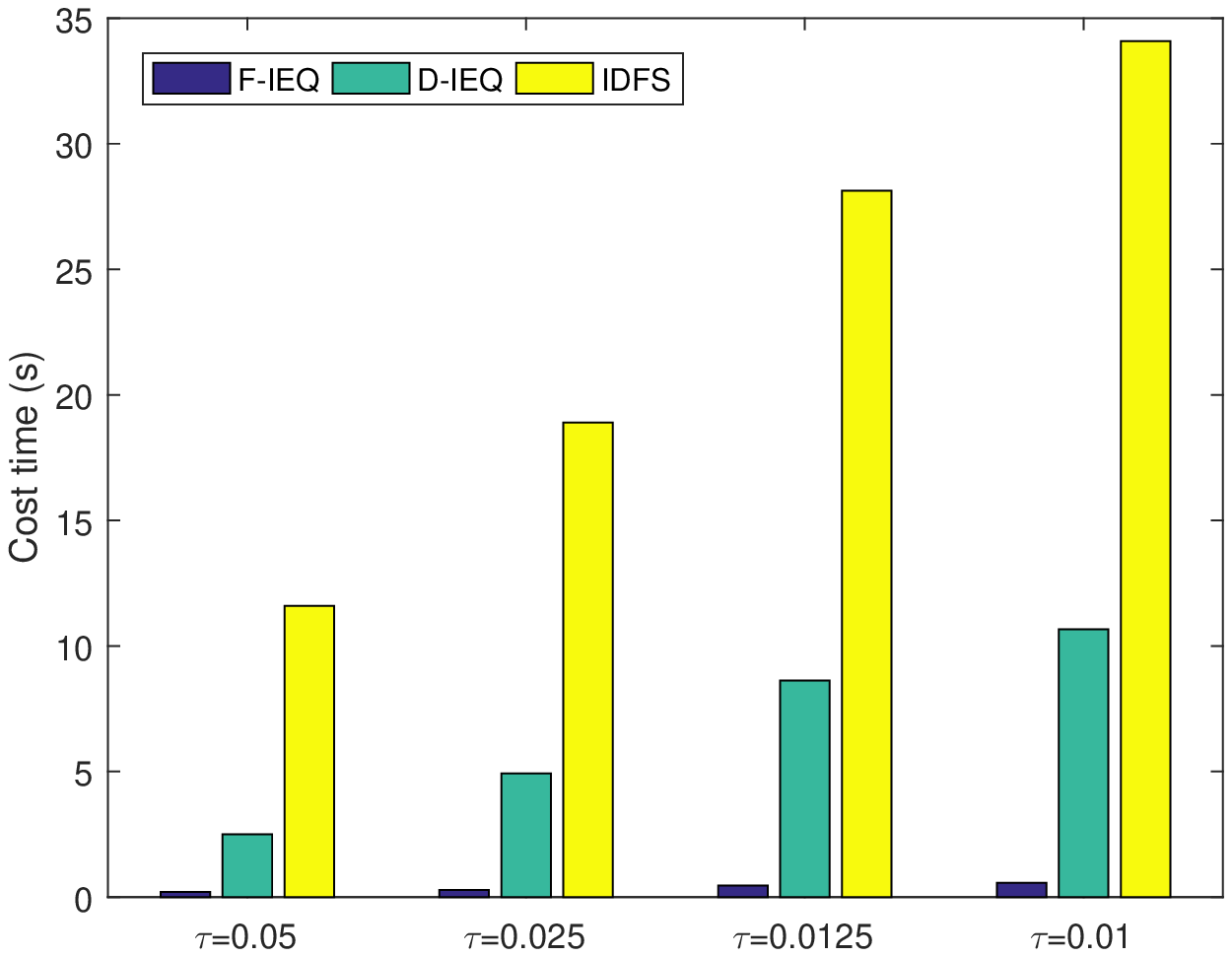}\\
{\footnotesize  \centerline {(a) $\alpha=1.3$}}
\end{minipage}
\begin{minipage}[t]{70mm}
\includegraphics[width=70mm]{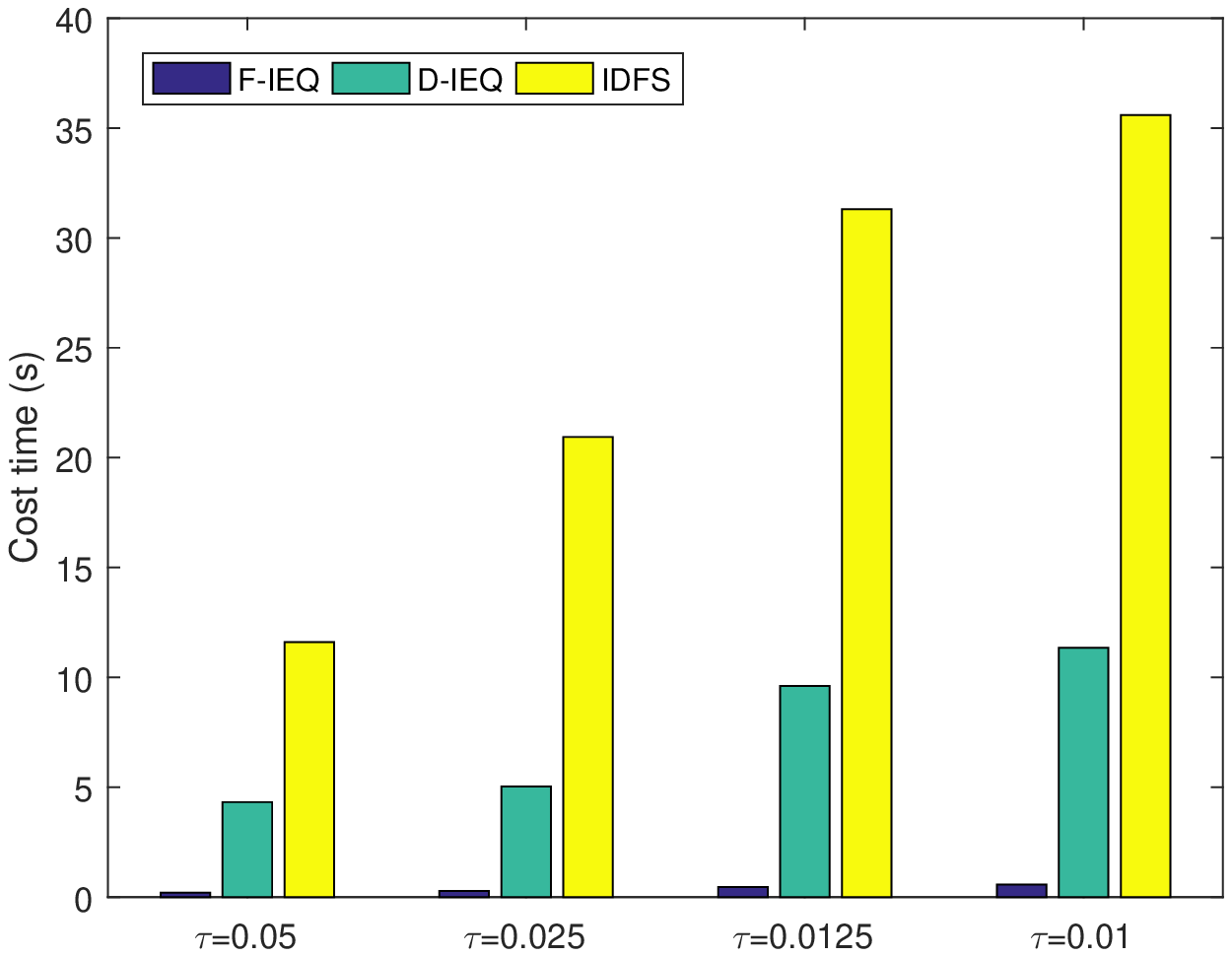}\\
{\footnotesize  \centerline {(b) $\alpha=2$}}
\end{minipage}
\caption{\small {\small{CPU time for the soliton with
different time steps till $T = 10$ under $h = 0.1$.} }}\label{fig522}
\end{figure}

\begin{table}[H]
\tabcolsep=8pt
\small
\renewcommand\arraystretch{1.2}
\centering
\caption{{The numerical errors and convergence orders of the IEQ-CN scheme at $T=1$.}}\label{Tab_2SG:2}
\begin{tabularx}{1.0\textwidth}{XXXXXX}\hline
{$\alpha$\ \ }&{$(h,\tau)$}  & {error} & {order} \\     
\hline

 {1.3} &{$(\frac{1}{5}$, $\frac{1}{50}$)}    & {1.5583e-03} & {-} \\
 {}    &{$(\frac{1}{10}$, $\frac{1}{100}$)}  & {3.8978e-04} & {1.9993} \\
 {}    &{$(\frac{1}{20}$, $\frac{1}{200}$)}  & {9.7441e-05} & {2.0000}  \\   
 {}    &{$(\frac{1}{40}$, $\frac{1}{400}$)}  & {2.4357e-05} & {2.0001} \\\hline
 {1.75}&{$(\frac{1}{5}$, $\frac{1}{50}$)}    & {2.4035e-03} & {-} \\
 {}    &{$(\frac{1}{10}$, $\frac{1}{100}$)}  & {5.9925e-04} & {2.0033} \\
 {}    &{$(\frac{1}{20}$, $\frac{1}{200}$)}  & {1.4969e-04} & {2.0011}  \\   
 {}    &{$(\frac{1}{40}$, $\frac{1}{400}$)}  & {3.7413e-05} & {2.0003} \\\hline
 {1.99}&{$(\frac{1}{5}$, $\frac{1}{50}$)}    & {2.7569e-03} & {-} \\
 {}    &{$(\frac{1}{10}$, $\frac{1}{100}$)}  & {6.8571e-04} & {2.0074} \\
 {}    &{$(\frac{1}{20}$, $\frac{1}{200}$)}  & {1.7119e-04} & {2.0019}  \\   
 {}    &{$(\frac{1}{40}$, $\frac{1}{400}$)}  & {4.2781e-05} & {2.0005} \\\hline
 {2}   &{$(\frac{1}{5}$, $\frac{1}{50}$)}    & {2.7689e-03} & {-} \\
 {}    &{$(\frac{1}{10}$, $\frac{1}{100}$)}  & {6.8864e-04} & {2.0075} \\
 {}    &{$(\frac{1}{20}$, $\frac{1}{200}$)}  & {1.7192e-04} & {2.0020}  \\   
 {}    &{$(\frac{1}{40}$, $\frac{1}{400}$)}  & {4.2963e-05} & {2.0006} \\\hline
\end{tabularx}
\end{table}

Second, we enlarge the computational domain $\Omega=(-40,40)$ and verify the discrete energy conservation law of the fully-discrete scheme. We take $h=0.1,$~$\tau=0.05$ and compute the discrete energy. Fig.2 shows the relative errors of energy $E$ for different values of fractional order $\alpha$.
The pictures demonstrate that the IEQ-CN  scheme preserves the energy very well in discrete sense.

\begin{figure}[H]
\centering\begin{minipage}[t]{70mm}
\includegraphics[width=60mm]{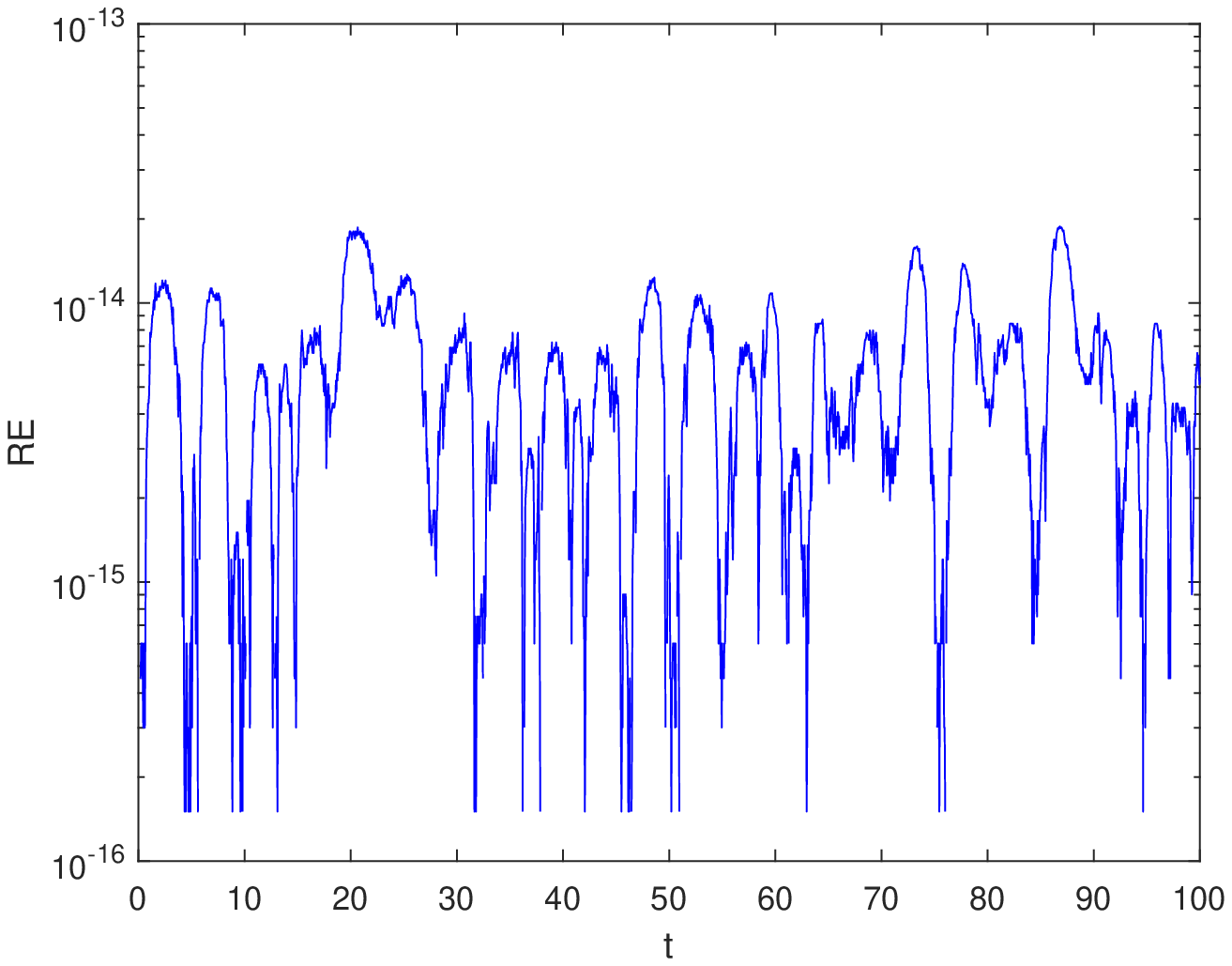}\\
{\footnotesize  \centerline {(a) $\alpha=1.3$}}
\end{minipage}
\begin{minipage}[t]{70mm}
\includegraphics[width=60mm]{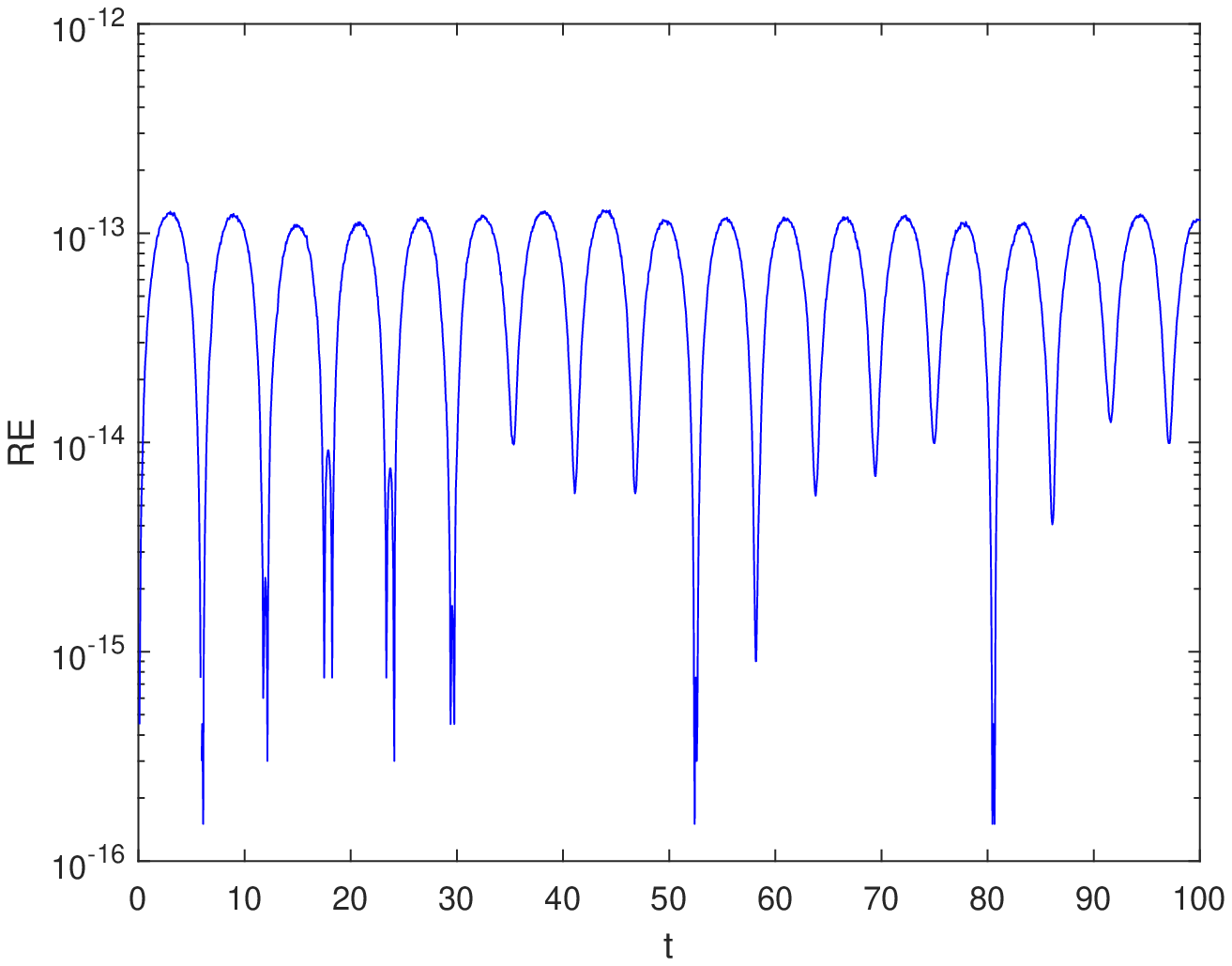}\\
{\footnotesize  \centerline {(b) $\alpha=1.75$}}
\end{minipage}
\begin{minipage}[t]{70mm}
\includegraphics[width=60mm]{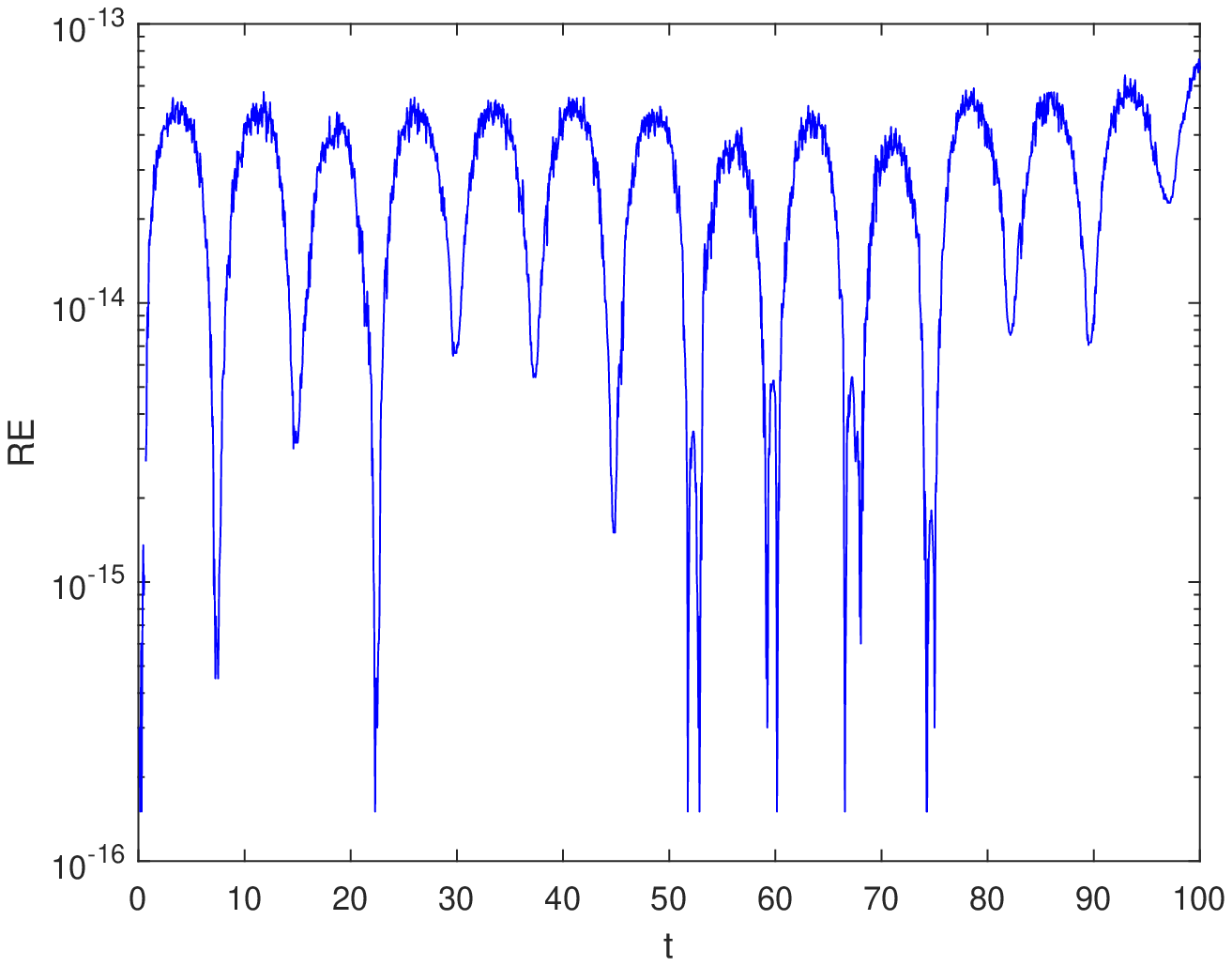}\\
{\footnotesize  \centerline {(c) $\alpha=1.99$}}
\end{minipage}
\begin{minipage}[t]{70mm}
\includegraphics[width=60mm]{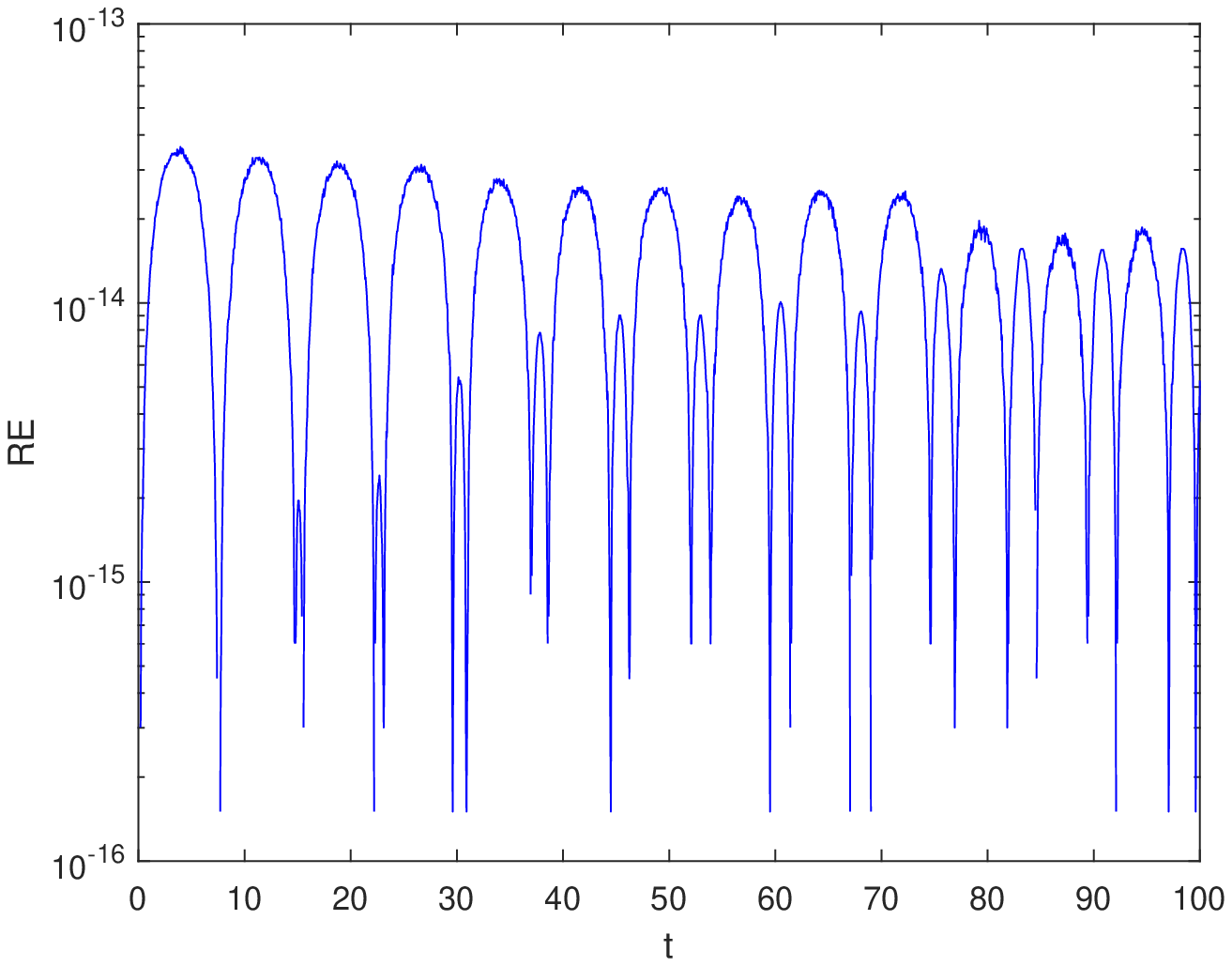}\\
{\footnotesize   \centerline {(d) $\alpha=2$}}
\end{minipage}
\caption{\small {The relative energy error with different $\alpha$ when $\omega=1.1$.}}\label{fig522}
\end{figure}

Last but not least, we select $\omega=1$ and pay attention to study the relationship between the evolution of the soliton and the fractional order $\alpha$  for the original system (\ref{FSG:eq:1.1}). Here we take the computation domain $\Omega=(-100,100)$. Without loss of generality, we choose $\alpha=1.1, 1.75, 1.99, 2$, $h=0.1$, $\tau=0.05$ and the numerical results are presented in Fig.3. Obviously, we can deduce that the shape of the soliton changes dramatically when
the fractional order $\alpha$ changes from 2.0 to 1.99. When $1< \alpha < 2$, the bigger the fractional order $\alpha$ is, the bigger the period of the soliton is.

\begin{figure}[H]
\centering\begin{minipage}[t]{70mm}
\includegraphics[width=60mm]{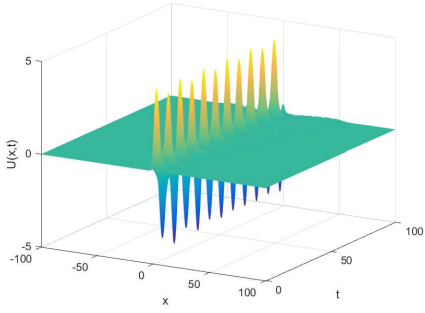}\\
{\footnotesize  \centerline {(a) $\alpha=1.1$}}
\end{minipage}
\begin{minipage}[t]{70mm}
\includegraphics[width=60mm]{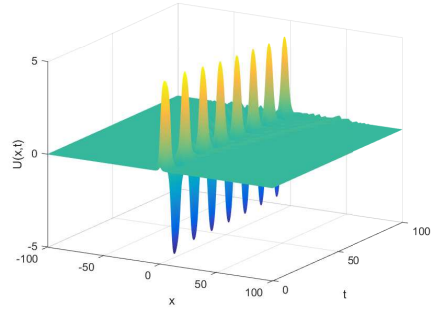}\\
{\footnotesize  \centerline {(b) $\alpha=1.75$}}
\end{minipage}
\begin{minipage}[t]{70mm}
\includegraphics[width=60mm]{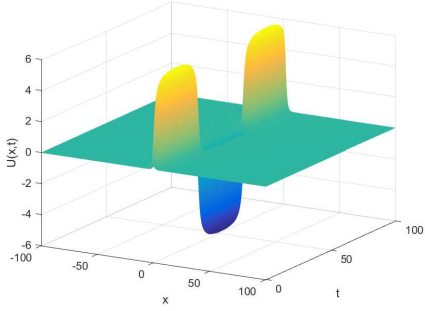}\\
{\footnotesize  \centerline {(c) $\alpha=1.99$}}
\end{minipage}
\begin{minipage}[t]{70mm}
\includegraphics[width=60mm]{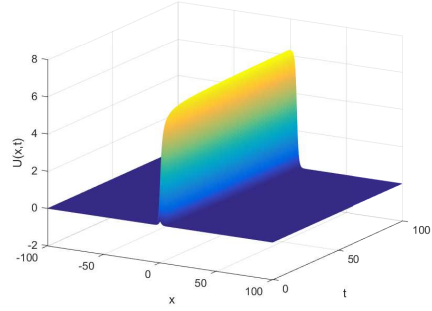}\\
{\footnotesize   \centerline {(d) $\alpha=2$}}
\end{minipage}
\caption{\small {\small{Evolution of the solitons with $\omega=1$ for different order $\alpha$.} }}\label{fig522}
\end{figure}

{\bf Example 5.2.} We consider the FSG equation with the initial conditions
\begin{align}
&\varphi(x)=3.2\text{sech}(x),\\
&\psi(x)=0.
\end{align}

In our computation, we take the Dirichlet boundary value
\begin{align*}
u(a,t)=u(b,t)=0,~~0\leq t \leq T.
\end{align*}

First, we set the space interval $\Omega=(-20,20)$ and test the accuracy of the fully-discrete scheme when $\alpha=1.3,~1.6,~1.9,~2$. As illustrated in Table. 2, the IEQ-CN scheme is second order of convergence in both time and space direction, which confirms the theoretical analysis.

\begin{table}[H]
\tabcolsep=8pt
\small
\renewcommand\arraystretch{1.2}
\centering
\caption{{The numerical errors and convergence orders of the IEQ-CN scheme at $T=1$.}}\label{Tab_2SG:2}
\begin{tabularx}{1.0\textwidth}{XXXXXX}\hline
{$\alpha$\ \ }&{$(h,\tau)$}  & {error} & {order} \\     
\hline

 {1.3} &{$(\frac{1}{5}$, $\frac{1}{50}$)}    & {4.3475e-03} & {-} \\
 {}    &{$(\frac{1}{10}$, $\frac{1}{100}$)}  & {1.0849e-03} & {2.0026} \\
 {}    &{$(\frac{1}{20}$, $\frac{1}{200}$)}  & {2.7117e-04} & {2.0003}  \\   
 {}    &{$(\frac{1}{40}$, $\frac{1}{400}$)}  & {6.7796e-05} & {1.9999} \\\hline
 {1.6}&{$(\frac{1}{5}$, $\frac{1}{50}$)}     & {5.1079e-03} & {-} \\
 {}    &{$(\frac{1}{10}$, $\frac{1}{100}$)}  & {1.2689e-03} & {2.0092} \\
 {}    &{$(\frac{1}{20}$, $\frac{1}{200}$)}  & {3.1678e-04} & {2.0020}  \\   
 {}    &{$(\frac{1}{40}$, $\frac{1}{400}$)}  & {7.9175e-05} & {2.0004} \\\hline
 {1.9}&{$(\frac{1}{5}$, $\frac{1}{50}$)}     & {5.1156e-03} & {-} \\
 {}    &{$(\frac{1}{10}$, $\frac{1}{100}$)}  & {1.2667e-03} & {2.0138} \\
 {}    &{$(\frac{1}{20}$, $\frac{1}{200}$)}  & {3.1601e-04} & {2.0031}  \\   
 {}    &{$(\frac{1}{40}$, $\frac{1}{400}$)}  & {7.8969e-05} & {2.0006} \\\hline
 {2}   &{$(\frac{1}{5}$, $\frac{1}{50}$)}    & {4.9566e-03} & {-} \\
 {}    &{$(\frac{1}{10}$, $\frac{1}{100}$)}  & {1.2273e-03} & {2.0139} \\
 {}    &{$(\frac{1}{20}$, $\frac{1}{200}$)}  & {3.0617e-04} & {2.0031}  \\   
 {}    &{$(\frac{1}{40}$, $\frac{1}{400}$)}  & {7.6510e-05} & {2.0006} \\\hline
\end{tabularx}
\end{table}

Second, we enlarge the computational domain $\Omega=(-40,40)$ and test the discrete energy conservation law of the IEQ-CN  scheme. The relative energy errors at different fractional order $\alpha$ = 1.3, 1.6, 1.9, 2
are presented in Fig. 4, where the numerical results are obtained with $h=\tau=0.05$.
 One can observe that the IEQ-CN  scheme preserves the energy well in discrete sense.

\begin{figure}[H]
\centering\begin{minipage}[t]{70mm}
\includegraphics[width=60mm]{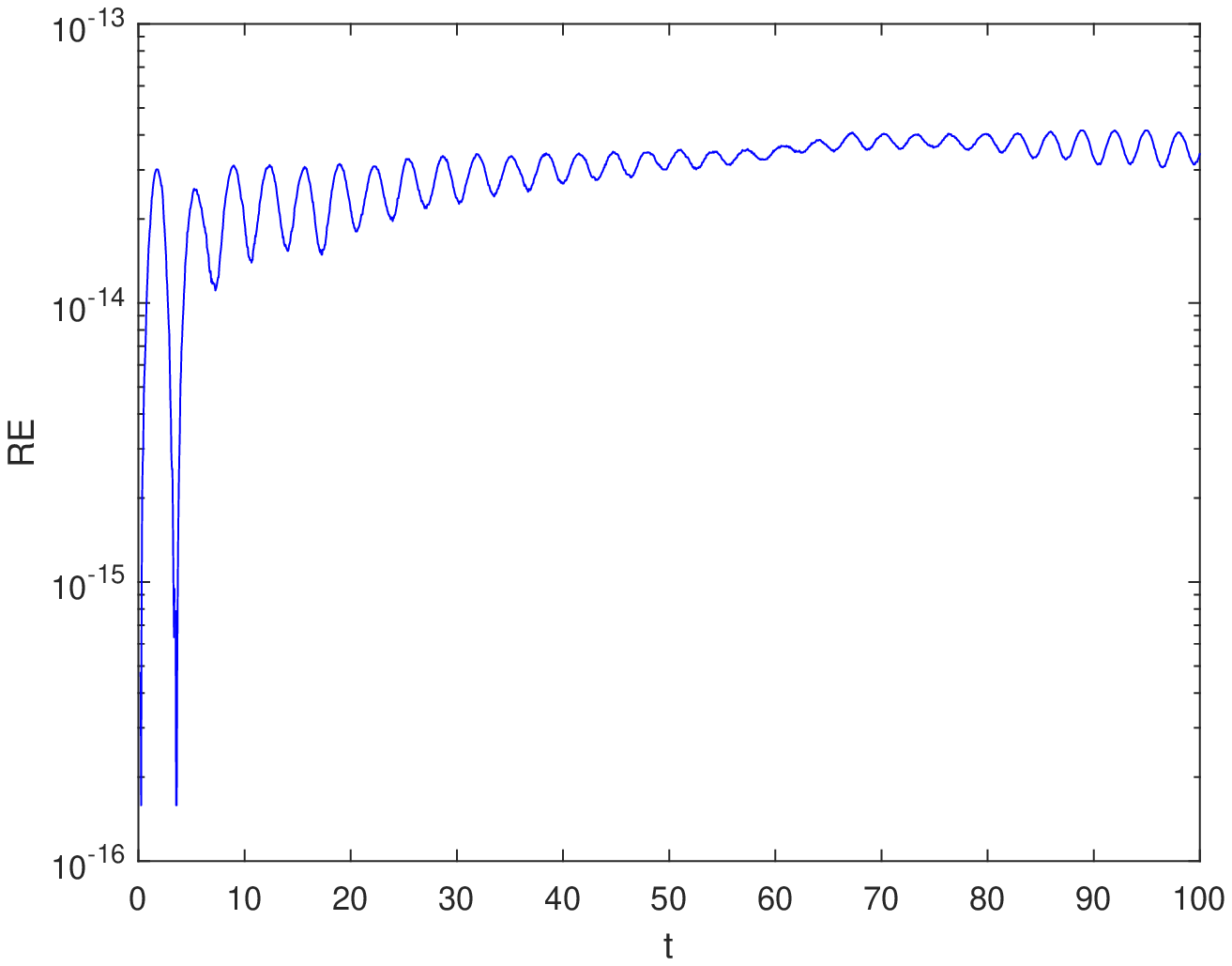}\\
{\footnotesize  \centerline {(a) $\alpha=1.3$}}
\end{minipage}
\begin{minipage}[t]{70mm}
\includegraphics[width=60mm]{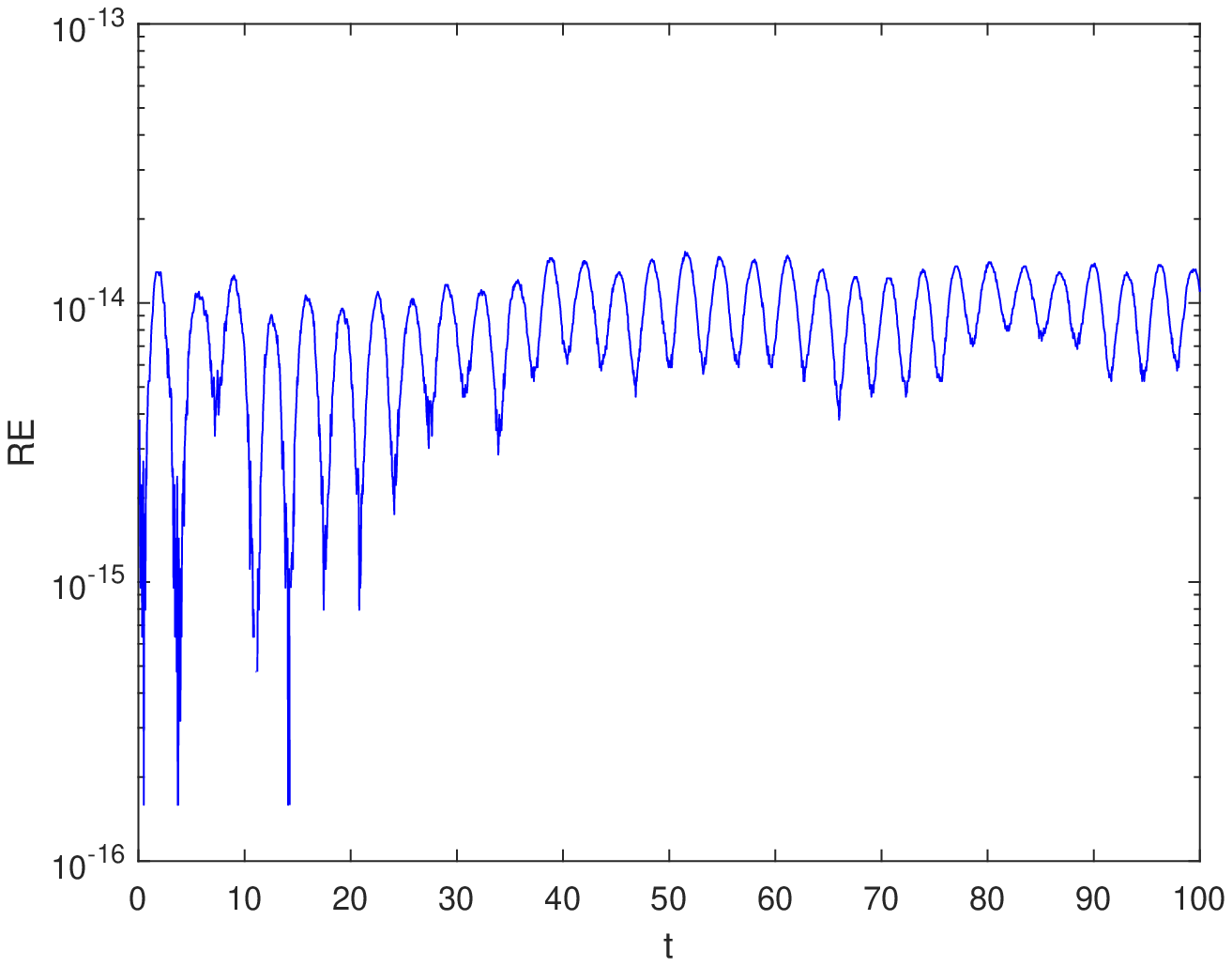}\\
{\footnotesize  \centerline {(b) $\alpha=1.6$}}
\end{minipage}
\begin{minipage}[t]{70mm}
\includegraphics[width=60mm]{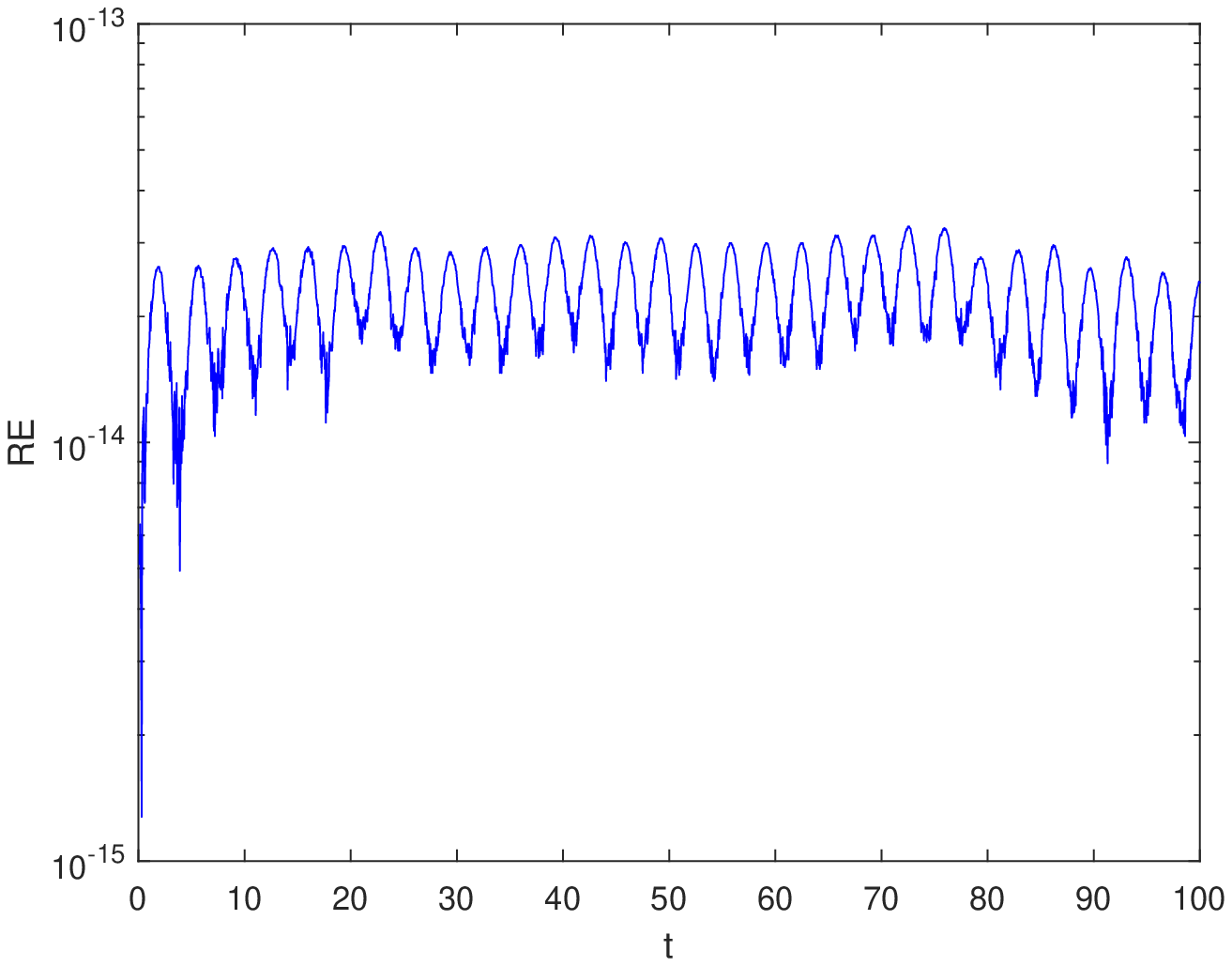}\\
{\footnotesize  \centerline {(c) $\alpha=1.9$}}
\end{minipage}
\begin{minipage}[t]{70mm}
\includegraphics[width=60mm]{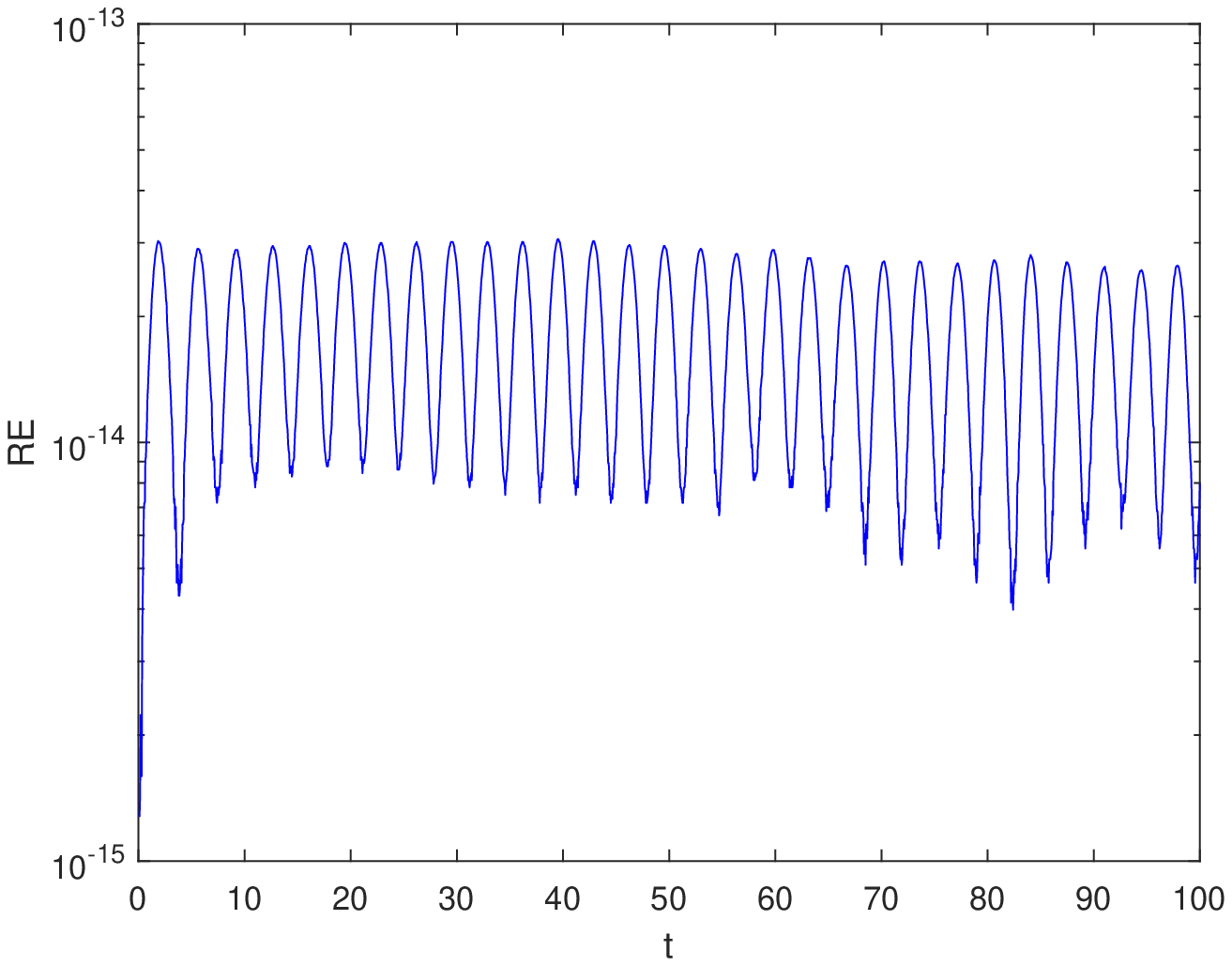}\\
{\footnotesize   \centerline {(d) $\alpha=2$}}
\end{minipage}
\caption{\small {The relative energy error with $h=\tau=0.05$ for different order $\alpha$.}}\label{fig522}
\end{figure}

Finally, we investigate the relationship between the fractional order $\alpha$ and the shape of the soliton for the problem with different $\alpha$. Here we take the computation domain $\Omega=(-100,100)$. The numerical solutions obtained by the IEQ-CN scheme with $h=0.1,~\tau=0.05$ are presented in Fig. 5. The results demonstrate that fractional order $\alpha$ will affect the shape of the soliton, and the shape of the soliton will change more quickly when $\alpha$ becomes smaller.

\vskip 1mm

\begin{figure}[H]
\centering\begin{minipage}[t]{70mm}
\includegraphics[width=60mm]{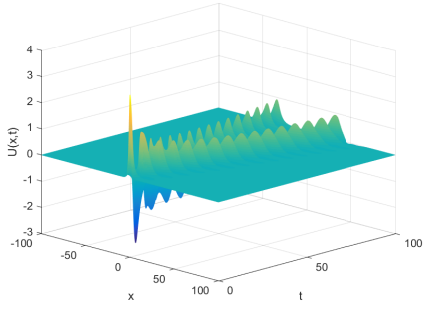}\\
{\footnotesize  \centerline {(a) $\alpha=1.3$}}
\end{minipage}
\begin{minipage}[t]{70mm}
\includegraphics[width=60mm]{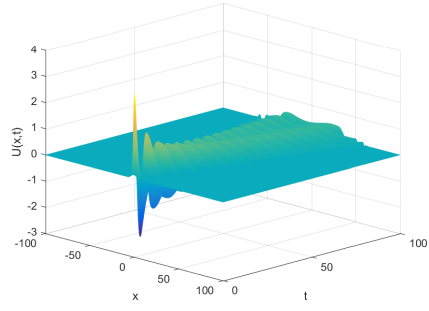}\\
{\footnotesize  \centerline {(b) $\alpha=1.6$}}
\end{minipage}
\begin{minipage}[t]{70mm}
\includegraphics[width=60mm]{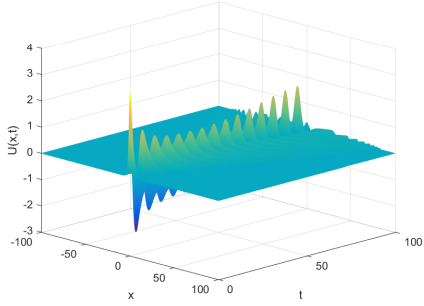}\\
{\footnotesize  \centerline {(c) $\alpha=1.9$}}
\end{minipage}
\begin{minipage}[t]{70mm}
\includegraphics[width=60mm]{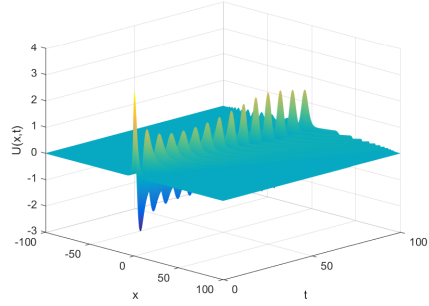}\\
{\footnotesize   \centerline {(d) $\alpha=2$}}
\end{minipage}
\caption{\small {\small{Evolution of the solitons for different order $\alpha$. }}}\label{fig522}
\end{figure}

\section{Conclusions}

In this paper, we derive the Hamiltonian formulation of the fractional sine-Gordon equation, and then construct a new difference scheme for the equation based on the invariant energy quadratization approach. Specifically, the scheme is linear, and can preserve discrete energy.
Theoretical analysis and numerical experiments indicate that the new scheme is efficient and accurate, and has desirable
energy conservation property. In addition, the proposed energy-preserving scheme can be generalized to other fractional equations, such as the nonlinear fractional Schr\"{o}dinger equation, the fractional Klein-Gordon-Schr\"{o}dinger equation, etc.

Recently, a new method which is termed as scalar auxiliary variable approach has been developed by Shen et al. for solving gradient flows \cite{p45,p46}. It
inherits all advantages of invariant energy quadratization approach
but also overcomes most of its shortcomings. Future research should be devoted to establishing the linear implicit energy-preserving scheme based on the scalar auxiliary variable approach
for fractional differential equations.

\section*{Acknowledgments}

This work is supported by the Postgraduate Research $\&$ Practice
Innovation Program of Jiangsu Province (Grant Nos. KYCX19\_0776), the National Natural Science Foundation of China (Grant No. 11771213, 61872422),
the National Key Research and Development Project of China (Grant No. 2016YFC0600310, 2018YFC0603500, 2018YFC1504205),
the Major Projects of Natural Sciences of University in Jiangsu Province of China (Grant No. 18KJA110003), and the Priority Academic Program Development of Jiangsu Higher Education Institutions.

\end{document}